\documentclass[11pt]{article}
\usepackage{amsmath,amsthm,amscd,amssymb, color, amsfonts, latexsym, geometry,textcomp, authblk, hyperref, hhline}
\usepackage{fullpage}
\usepackage{bm}
\usepackage{latexsym}
\usepackage{graphicx}
\usepackage{tikz}
\usepackage{url}
\usepackage{tikz-cd}
\usetikzlibrary{arrows, matrix}
\usepackage[all]{xy}
\usepackage{multirow}
\usepackage{mathtools}

\usepackage{tocloft}

\newlength\mylen

\renewcommand\cftpartpresnum{Part~}
\newcommand{\PSLt}{{\mathrm{PSL}_2}}

\providecommand{\keywords}[1]{\textbf{\textit{Keywords: }} #1}

\providecommand{\msc}[1]{\textbf{\textit{2010 Mathematics Subject Classification: }} #1}

\newtheorem{problem}{Problem}

\newcommand{\QQ}{{\mathbb{Q}}}
\newcommand{\ZZ}{{\mathbb{Z}}}
\newcommand{\RR}{{\mathbb{R}}}
\newcommand{\FF}{{\mathbb{F}}}
\newcommand{\A}{\mathbb{A}}
\renewcommand{\O}{\mathcal{O}}

\renewcommand{\char}{\text{char}}
\newcommand{\Cl}{{\mathrm{Cl}}}
\newcommand{\tr}{{\mathrm{tr}}}

\newcommand{\Cay}{\mathrm{Cay}}
\newcommand{\Ql}{{\QQ_l}}
\newcommand{\Zl}{{\ZZ_l}}

\newcommand{\PGL}{{\mathrm{PGL}}}
\newcommand{\PGLt}{{\mathrm{PGL}_2}}

\newcommand{\tbtmtx}[4]{\left(\begin{array}{cc}
#1 & #2\\
#3 & #4
\end{array}\right)}
\newcommand{\Span}{{\mathrm{Span}}}
\newcommand{\Aut}{{\mathrm{Aut}}}

\DeclarePairedDelimiter\floor{\lfloor}{\rfloor}

\newtheorem{lemma}{Lemma}[section]
\newtheorem{proposition}[lemma]{Proposition}
\newtheorem{theorem}[lemma]{Theorem}
\newtheorem{definition}[lemma]{Definition}
\newtheorem{corollary}[lemma]{Corollary}

\theoremstyle{definition}
\newtheorem{example}[lemma]{Example}

\newtheorem{remark}[lemma]{Remark}

\date{}

\title{Ramanujan graphs in cryptography}
\author[1]{Anamaria Costache}
\author[2]{Brooke Feigon}
\author[3]{Kristin Lauter}
\author[4]{Maike Massierer}
\author[5]{Anna Pusk\'{a}s}
\affil[1]{\footnotesize Department of Computer Science, University of Bristol, Bristol, UK, \href{mailto:anamaria.costache@bristol.ac.uk}{anamaria.costache@bristol.ac.uk}}
\affil[2]{\footnotesize Department of Mathematics, The City College of New York, CUNY, NAC 8/133, New York, NY 10031,  \href{mailto:bfeigon@ccny.cuny.edu}{bfeigon@ccny.cuny.edu} \thanks{Partially supported by National Security Agency grant H98230-16-1-0017  and PSC-CUNY.}}
\affil[3]{\footnotesize Microsoft Research, One Microsoft Way, Redmond, WA 98052, \href{mailto:klauter@microsoft.com}{klauter@microsoft.com}}
\affil[4]{\footnotesize School of Mathematics and Statistics, University of New South Wales, Sydney NSW  2052, Australia, \href{mailto:maike.massierer@gmail.com}{maike.massierer@gmail.com} \thanks{Partially supported by Australian Research Council grant DP150101689.}}
\affil[5]{\footnotesize Department of Mathematics \& Statistics, University of Massachusetts, Amherst, MA 01003, \href{mailto:puskas@math.umass.edu}{puskas@math.umass.edu} }

\begin{document}

\maketitle

\begin{abstract} In this paper we study the security of a proposal for Post-Quantum Cryptography  from both a number theoretic and cryptographic perspective. Charles--Goren--Lauter in 2006~\cite{lauter:2009}  proposed two hash functions based on the hardness of finding paths in Ramanujan graphs. One is  based on Lubotzky--Phillips--Sarnak (LPS) graphs and the other one is based on Supersingular Isogeny Graphs.  A 2008 paper by Petit--Lauter--Quisquater breaks the hash function based on LPS graphs. On the Supersingular Isogeny Graphs proposal, recent work has continued to build cryptographic applications on the hardness of finding isogenies between supersingular elliptic curves. A 2011 paper by De Feo--Jao--Pl\^ut proposed a cryptographic system based on Supersingular Isogeny Diffie--Hellman as well as a set of five hard problems.  In this paper we show that the security of the SIDH proposal relies on the hardness of the SSIG path-finding problem introduced in~\cite{lauter:2009}. In addition, similarities between the number theoretic ingredients in the LPS and Pizer constructions suggest that the hardness of the path-finding problem in the two graphs may be linked. By viewing both graphs from a number theoretic perspective, we identify the similarities and differences between the Pizer and LPS graphs.
\end{abstract}


\msc{Primary: 05C25, 14G50; Secondary: 22F70, 11R52}

\section{Introduction}
Supersingular Isogeny Graphs were proposed for use in cryptography in 2006 by Charles, Goren, and Lauter~\cite{lauter:2009}.  Supersingular isogeny graphs are examples of Ramanujan graphs, i.e.\ optimal expander graphs.  This means that relatively {\it short} walks on the graph approximate the uniform distribution, i.e. walks of length approximately equal to the logarithm of the graph size. Walks on expander graphs are often used as a good source of randomness in computer science, and the reason for using {\it Ramanujan} graphs is to keep the path length short.  But the reason these graphs are important for cryptography is that {\it finding paths} in these graphs, i.e.\ {\it routing,} is hard: there are no known subexponential algorithms to solve this problem, either classically or on a quantum computer.  For this reason, systems based on the hardness of problems on Supersingular Isogeny Graphs are currently under consideration for standardization in the NIST Post-Quantum Cryptography (PQC) Competition~\cite{PQC}.

\cite{lauter:2009} proposed a general construction for cryptographic hash functions based on the hardness of inverting a walk on a graph.  The path-finding problem is the following: given fixed starting and ending vertices representing the start and end points of a walk on the graph of a fixed length, find a path between them.  A hash function can be defined by using the input to the function as directions for walking around the graph: the output is the label for the ending vertex of the walk. Finding collisions for the hash function is equivalent to finding cycles in the graph, and finding pre-images is equivalent to path-finding in the graph. Backtracking is not allowed in the walks by definition, to avoid trivial collisions.

In~\cite{lauter:2009}, two concrete examples of families of optimal expander graphs (Ramanujan graphs) were proposed, the so-called Lubotzky--Phillips--Sarnak (LPS) graphs~\cite{lubotzky:1988}, and the Supersingular Isogeny Graphs (Pizer) \cite{pizer:1998}, where the path finding problem was supposed to be hard.   Both graphs were proposed and presented at the 2005 and 2006 NIST Hash Function workshops, but the LPS hash function was quickly attacked and broken in two papers in 2008, a collision attack~\cite{Zemor-Tillich} and a pre-image attack~\cite{SCN:PetLauQui08}.  The preimage attack gives an algorithm to efficiently find paths in LPS graphs, a problem which had been open for several decades. The PLQ path-finding algorithm uses the explicit description of the graph as a Cayley graph in $\PSLt(\FF_p)$, where vertices are $2 \times 2$ matrices with entries in $\FF_p$ satisfying certain properties.  Given the swift discovery of attacks on the LPS path-finding problem, it is natural to investigate whether this approach is relevant to the path-finding problem in Supersingular Isogeny (Pizer) Graphs.

In 2011, De Feo--Jao--Pl\^ut \cite{defeo:2014}  devised a cryptographic system based on supersingular isogeny graphs, proposing a Diffie--Hellman protocol as well as a set of five hard problems related to the security of the protocol.  It is natural to ask what is the relation between the problems stated in~\cite{defeo:2014}  and the path-finding problem on Supersingular Isogeny Graphs proposed in~\cite{lauter:2009}.

In this paper we explore these two questions related to the security of cryptosystems based on these Ramanujan graphs. In Part \ref{part1} of the paper, we study the relation between the hard problems proposed by De Feo--Jao--Pl\^ut and the hardness of the Supersingular Isogeny Graph problem which is the foundation for the CGL hash function. In Part \ref{part2} of the paper, we study the relation between the Pizer and LPS graphs by viewing both from a number theoretic perspective.

In particular, in Part \ref{part1} of the paper, we clearly explain how the security of the Key Exchange protocol relies on the hardness of the path-finding problem in SSIG, proving a reduction (Theorem ~\ref{lemma::reduction})  between the Supersingular Isogeny Diffie Hellmann (SIDH) Problem and the path-finding problem in SSIG.   Although this fact and this theorem may be clear to the experts (see for example the comment in the introduction to a recent paper on this topic~\cite{Menezes}), this reduction between the hard problems is not written anywhere in the literature.  Furthermore, the Key Exchange (SIDH) paper~\cite{defeo:2014}  states 5 hard problems, including (SSCDH), with relations proved between some but not all of them, and mentions the paper \cite{lauter:2009} only in passing (on page 17), with no clear statement of the relationship to the overarching hard problem of path-finding in SSIG.

Our Theorem~\ref{lemma::reduction} clearly shows the fact that the security of the proposed post-quantum key exchange relies on the hardness of the path-finding problem in SSIG stated in~\cite{lauter:2009}. Theorem~\ref{thm::number-chains} counts the chains of isogenies of fixed length. Its proof relies on elementary group theory results and facts about isogenies, proved in Section \ref{sec::correspondence}. 

In Part \ref{part2} of the paper, we examine the LPS and Pizer graphs from a number theoretic perspective with the aim of highlighting the similarities and differences between the constructions.

Both the LPS and Pizer graphs considered in~\cite{lauter:2009} can be thought of as graphs on
\begin{equation}\label{eq:intro_double_coset}
\Gamma \backslash \PGLt(\Ql)/\PGLt(\Zl),
\end{equation}
where $\Gamma$ is a discrete cocompact subgroup, where $\Gamma$ is obtained from a quaternion algebra $B.$ We show how different input choices for the construction lead to different graphs. In the LPS construction one may vary $\Gamma$ to get an infinite family of Ramanujan graphs. In the Pizer construction one may vary $B$ to get an infinite family. In the LPS case, we always work in the Hamiltonian quaternion algebra. For this particular choice of algebra we can 
rewrite the graph
as a Cayley graph. This explicit description is key for breaking the LPS hash function.  For the Pizer graphs we do not have such a description. On the Pizer side the graphs may, via Strong Approximation, be viewed as graphs on ad\`{e}lic double cosets which are in turn the class group of an order of $B$ that is related to the cocompact subgroup $\Gamma$. From here one obtains an isomorphism with supersingular isogeny graphs. For LPS graphs the local double cosets are also isomorphic to ad\`{e}lic double cosets, but  
in this case the corresponding set of ad\`{e}lic double cosets is smaller relative to the quaternion algebra  
and we do not have the same chain of isomorphisms.

Part \ref{part2} has the following outline. Section \ref{subsect:LPS} follows \cite{lubotzky:2010} and presents the construction of LPS graphs from three different perspectives: as a Cayley graph, in terms of local double cosets, and, to connect these two, as a quotient of an infinite tree. The edges of the LPS graph are explicit in both the Cayley and local double coset presentation. In Section \ref{app:sect:explicitisom} we give an explicit bijection between the natural parameterizations of the edges at a fixed vertex. Section \ref{subsect:StrongApprox} is about Strong Approximation, the main tool connecting the local and adelic double cosets for both LPS and Pizer graphs. Section \ref{sect:Pizer} follows \cite{pizer:1998} and summarizes Pizer's construction. The different input choices for LPS and Pizer constructions impose different restrictions on the parameters of the graph, such as the degree. $6$-regular graphs exist in both families. In Section \ref{subsect:Pizer_primes} we give a set of congruence conditions for the parameters of the Pizer construction that produce a $6$-regular graph. In Section \ref{subsect:LPS_Pizer_compare} we summarize the similarities and differences between the two constructions.

\subsection{Acknowledgments}
This project was initiated at the Women in Numbers 4 (WIN4) workshop at the Banff International Research Station in August, 2017. The authors would like to thank BIRS and the WIN4 organizers. In addition, the authors would like to thank the Clay Mathematics Institute, PIMS, Microsoft Research, the Number Theory Foundation and the NSF-HRD 1500481 - AWM ADVANCE  grant for supporting the workshop.  We thank John Voight, Scott Harper, and Steven Galbraith for helpful conversations, and the anonymous referees for many helpful suggestions and edits.

\part{\Large Cryptographic applications of supersingular isogeny graphs}\label{part1}
In this section we investigate the security of the \cite{defeo:2014}
key-exchange protocol.  We show a reduction to the path-finding problem in supersingular isogeny
graphs stated in \cite{lauter:2009}.  The hardness of this problem is the
basis for the CGL cryptographic hash function, and we show here that if this problem is not hard, then the
key exchange presented in \cite{defeo:2014} is not secure.

We begin by recalling some basic facts about
isogenies of elliptic curves and the key-exchange construction. Then, we give a reduction between two hardness assumptions.  This reduction is based on a correspondence between a path representing the composition of $m$ isogenies of degree $\ell$ and an isogeny of degree
$\ell^m$.

\section{Preliminaries} \label{sec::crypto-prelim}

We start by recalling some basic and well-known results about isogenies. They
can all be found in \cite{silverman}. We try to be as concrete and constructive
as possible, since we would like to use these facts to do computations.

An elliptic curve is a curve of genus one with a specific base point $\O$. This
latter can be used to define a group law. We will not go into the details of
this, see for example \cite{silverman}. If $E$ is an elliptic curve defined over a field
$K$ and $\char(\bar{K})\neq 2,3$, we can write the equation of $E$ as

\[ E: y^2 = x^3 +a\cdot x + b, \]
where $a, b \in K$.
Two important quantities related to an elliptic curve are its discriminant
$\Delta$ and its $j$-invariant, denoted by $j$. They are defined as follows.

\[ \Delta = 16 \cdot (4\cdot a^3 + 27\cdot b^2) \quad\text{and}\quad j = -1728\cdot\frac{a^3}{\Delta}.  \]
Two elliptic curves are isomorphic over $\bar{K}$ if and only if they have the same $j$-invariant.

\begin{definition} Let $E_0$ and $E_1$ be two elliptic curves. An isogeny from
$E_0$ to $E_1$ is a surjective morphism
\[ \phi: E_0 \rightarrow E_1,  \]
which is a group homomorphism.
\end{definition}
An example of an isogeny is the multiplication-by-$m$ map $[m]$,
\begin{align*}
[m] : E &\rightarrow E \\
      P & \mapsto m\cdot P.
\end{align*}

The degree of an isogeny is defined as the degree of the finite extension
$\bar{K}(E_0)/\phi^*(\bar{K}(E_1))$, where $\bar{K}(*)$ is the
function field of the curve, and $\phi^*$ is the map of function fields induced by the isogeny
$\phi$. By convention, we set
\[ \deg([0]) = 0. \]
The degree map is multiplicative under composition of isogenies:
\[ \deg(\phi \circ \psi) = \deg(\phi) \cdot \deg(\psi)  \]
for all chains $E_0 \xrightarrow{\phi} E_1 \xrightarrow{\psi} E_2$, and for an integer $m > 0$, the multiplication-by-$m$ map has degree $m^2$.

\begin{theorem}{\cite{silverman}} Let $E_0 \rightarrow E_1$ be an
	isogeny of degree $m$. Then, there exists a unique isogeny
	\[ \hat{\phi}: E_1 \rightarrow E_0   \]
	such that $\hat{\phi} \circ \phi = [m]$ on $E_0$, and $\phi \circ
	\hat{\phi}=[m]$ on $E_1$. We call $\hat{\phi}$ the dual isogeny to $\phi$.
	We also have that
	\[ \deg(\hat{\phi}) = \deg(\phi).  \]
\end{theorem}

For an isogeny $\phi$, we say $\phi$ is separable if the field extension
$\bar{K}(E_0)/\phi^*(\bar{K}(E_1))$ is separable. We then have the following lemma.
\begin{lemma} Let $\phi: E_0 \rightarrow E_1$ be a separable isogeny. Then
	\[ \deg(\phi) = \#\ker(\phi).   \]
\end{lemma}
In this paper, we only consider separable isogenies and frequently use this convenient fact.
From the above, it follows that a point $P$ of order $m$ defines an isogeny $\phi$
of degree $m$,
\[ \phi: E \rightarrow E/\langle P \rangle.  \]
We will refer to such an isogeny as a cyclic isogeny (meaning that its kernel is a cyclic subgroup of $E$). For $\ell$ prime, we also say that
two curves $E_0$ and $E_1$ are $\ell$-isogenous if there exists an isogeny
$\phi: E_0 \rightarrow E_1$ of degree $\ell$.

We define $E[m]$,
the $m$-torsion subgroup of $E$, to be the kernel of the multiplication-by-$m$ map.
If $\char(K)>0$ and $m \geq 2$ is an integer coprime to $\char(K)$, or if
$\char(K)=0$, then the points of $E[m]$ are
\[ E[m] = \{ P \in E(\bar{K}): m \cdot P = \O  \} \cong \ZZ/m\ZZ \times \ZZ/m\ZZ.   \]

If an elliptic curve $E$ is defined over a field of characteristic $p>0$ and its
endomorphism ring over $\bar{K}$ is an order in a quaternion algebra, we say that $E$ is
supersingular.  Every isomorphism class over $\bar{K}$ of supersingular elliptic curves in characteristic $p$ has a representative defined over $\FF_{p^2}$, thus we will often let $K=\FF_{p^2}$ (for some fixed prime $p$).

We mentioned above that an $\ell$-torsion point $P$ induces an isogeny of degree $\ell$. More generally, a finite subgroup
$G$ of $E$ generates a unique isogeny of degree $\# G$, up to automorphism.

Supersingular isogeny graphs were introduced into cryptography in \cite{lauter:2009}. To define a supersingular isogeny graph, fix a finite field $K$ of characteristic $p$, a supersingular elliptic curve $E$ over $K$, and a prime $\ell \ne p$. Then the corresponding isogeny graph is constructed as follows. The vertices are the $\bar{K}$-isomorphism classes of elliptic curves which are $\bar{K}$-isogenous to $E$. Each vertex is labeled with the $j$-invariant of the curve. The edges of the graph correspond to the $\ell$-isogenies between the elliptic curves. As the vertices are isomorphism classes of elliptic curves, isogenies that differ by composition with an automorphism of the image are identified as edges of the graph. I.e. if $E_0,E_1$ are $\bar{K}$-isogenous elliptic curves, $\phi:E_0\rightarrow E_1$ is an $\ell$-isogeny and $\epsilon\in \Aut(E_1)$ is an automorphism, then $\phi$ and $\epsilon\circ \phi $ are identified and correspond to the same edge of the graph.  

If $p \equiv 1\mod {12}$, we can uniquely identify an isogeny with its dual to make it an undirected graph.  It is a multigraph in the sense that there can be multiple edges if no extra conditions are imposed on $p$. Three important properties of these graphs follow from deep theorems in number theory:
\begin{enumerate}
\item The graph is connected for any $\ell \ne p$
(special case of~\cite[Theorem 4.1]{lauter:2009r}).
\item A supersingular isogeny graph has roughly $p/12$ vertices. ~\cite[Theorem 4.1]{silverman}
\item Supersingular isogeny graphs are optimal expander graphs, in particular they are Ramanujan. (special case of~\cite[Theorem 4.2]{lauter:2009r}).
\end{enumerate}

\begin{remark} \label{Rem: multiple edges}
In order to avoid trivial collisions in cryptographic hash functions based on isogeny graphs, it is best if the graph has no short cycles. Charles, Goren, and Lauter show in \cite{lauter:2009} how to ensure that isogeny graphs do not have short cycles by carefully choosing the finite field one works over. For example, they compute that a $2$-isogeny graph does not have double edges (i.e.\ cycles of length $2$) when working over $\FF_p$ with $p \equiv 1 \bmod 420$. Similarly, we computed that a $3$-isogeny graph does not have double edges for $p \equiv 1 \bmod 9240$. Given that $420 = 2^2 \cdot 3 \cdot 5 \cdot 7$ and $9240 = 2^3 \cdot 3 \cdot 5 \cdot 7 \cdot 11$, we conclude that neither the 2-isogeny graph nor the 3-isogeny graph has double edges for $p \equiv 1 \bmod 9240$.

For our experiments (described in Section \ref{sec::correspondence}), we were interested in studying short walks, for example of length $4$, in a setting relevant to the Key-Exchange protocol described below.  The smallest prime $p$ with the property $p \equiv 1 \bmod 9240$ that also satisfies $2^4\cdot3^4 \mid p-1$ is
$$p = 2^4\cdot 3^4 \cdot 5 \cdot 7 \cdot 11 +1.$$
\end{remark}

\section{The \cite{defeo:2014} key-exchange}
Let $E$ be a supersingular elliptic curve defined over $\FF_{p^2}$,
where $p =\ell_A^n\cdot \ell_B^m\pm 1$, $\ell_A$ and $\ell_B$ are primes, and $n \approx m$ are approximately equal.
We have players $A$ (for Alice) and $B$ (for Bob), representing the two parties who wish to engage in a key-exchange protocol with the goal of establishing a shared secret key by communicating via a (possibly) insecure channel.
The two players $A$ and $B$ generate their public parameters by each picking two points $P_A$, $Q_A$ such that $\langle P_A, Q_A \rangle = E[\ell_A^n]$ (for $A$), and two points $P_B$, $Q_B$ such that $\langle P_B, Q_B \rangle = E[\ell_B^m]$ (for $B$).

Player $A$ then secretly picks two random integers $0\leq m_A, n_A< \ell_A^n$. These two integers (and the isogeny they generate) will be player $A$'s secret parameters. $A$ then computes the isogeny $\phi_A$
\[ E \xrightarrow{\phi_A} E_A :=E/\langle[m_A]P_A + [n_A]Q_A\rangle.   \]

Player $B$ proceeds in a similar fashion and secretly picks $0\leq m_B, n_B < \ell_B^m$.
Player $B$ then generates the (secret) isogeny
\[ E \xrightarrow{\phi_B} E_B :=E/\langle[m_B]P_B + [n_B]Q_B\rangle.   \]

So far, $A$ and $B$ have constructed the following diagram.
 
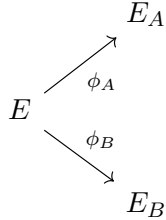
\[
 \begin{tikzcd}{}
 & E_A\arrow[leftarrow]{dl}{\phi_A}
 &   \\
 E \arrow{dr}{\phi_B}
 & & \\
 & E_B
 &
 \end{tikzcd}
 \]
 To complete the diamond, we proceed to the exchange part of the protocol.
 Player
 $A$ computes the points $\phi_A(P_B)$ and $\phi_A(Q_B)$ and sends $\{ \phi_A(P_B), \phi_A(Q_B), E_A \}$ along
 to player $B$. Similarly, player $B$ computes and sends $\{ \phi_B(P_A), \phi_B(Q_A), E_B \}$ to player $A$. Both players now have enough information to	construct the following diagram,

 \begin{equation}\label{diagram}
 \begin{tikzcd}{}
 & E_A\arrow[leftarrow]{dl}{\phi_A}\arrow[rightarrow]{dr}{\phi'_{A}}
 &   \\
 E \arrow{dr}{\phi_B}
 & & E_{AB} \arrow[leftarrow]{dl}{\phi'_{B}} \\
 & E_B
 &
 \end{tikzcd}
 \end{equation}
where
\[ E_{AB} \cong E/\langle [m_A]P_A + [n_A]Q_A, [m_B]P_B + [n_B]Q_B \rangle.   \]
Player $A$ can use the knowledge of the secret information $m_A$ and $n_A$ to compute the isogeny $\phi'_{B}$, by quotienting $E_B$ by $\langle[m_A]\phi_B(P_A) + [n_A]\phi_B(Q_A) \rangle$ to obtain $E_{AB}$.
Player $B$ can use the knowledge of the secret information $m_B$ and $n_B$ to compute the isogeny $\phi'_{A}$, by quotienting $E_A$ by $\langle[m_B]\phi_A(P_B) + [n_B]\phi_A(Q_B) \rangle$ to obtain $E_{AB}$.
A separable isogeny is determined by its kernel, and so both ways of going around the diagram from $E$ result in computing the same elliptic curve $E_{AB}$.

The players then use the $j$-invariant of the curve $E_{AB}$ as a shared secret.

\begin{remark} Given a list of points specifying a kernel, one can explicitly compute
	the associated isogeny using V\'{e}lu's formulas \cite{velu:1971}. In principle, this is how the
	two parties engaging in the key-exchange above can compute $\phi_A$, $\phi_B$, $\phi'_{A}$, $\phi'_{B}$
     \cite{velu:1971}.  However, in practice for cryptographic size subgroups, this would be impossible, and thus
    a different approach is taken, based on breaking the isogenies into $n$ (resp. $m$) steps, each of degree $\ell_A$ (resp. $\ell_B$).  This equivalence will be explained below.
\end{remark}

\subsection{Hardness assumptions}

The security of the key-exchange protocol is based on the following hardness assumption, which was introduced in \cite{defeo:2014} and called the Supersingular Computational Diffie--Hellman (SSCDH) problem.

\begin{problem}(Supersingular Computational Diffie--Hellman (SSCDH)):
\label{pbm:sscdh}
 Let $p$, $\ell_A$, $\ell_B$, $n$, $m$, $E$, $E_A$, $E_B$, $E_{AB}$, $P_A$, $Q_A$, $P_B$, $Q_B$ be as above.

Let $\phi_A$ be an isogeny from $E$ to $E_A$ whose kernel is equal to $\langle [m_A]P_A + [n_A]Q_A \rangle$, and let $\phi_B$ be an isogeny from $E$ to $E_B$ whose kernel is equal to $\langle [m_B]P_B + [n_B]Q_B \rangle$, where $m_A$,$n_A$ (respectively $m_B$,$n_B$) are integers chosen at random between $0$ and $\ell_A^m$ (respectively $\ell_B^n$), and not both divisible by $\ell_A$ (resp. $\ell_B$).

Given the curves $E_A$, $E_B$ and the points $\phi_A(P_B)$, $\phi_A(Q_B)$, $\phi_B(P_A)$, $\phi_B(Q_A)$, find the $j$-invariant of
\[ E_{AB} \cong E/\langle [m_A]P_A + [n_A]Q_A, [m_B]P_B + [n_B]Q_B \rangle;   \]
 see diagram \eqref{diagram}.
\end{problem}

In~\cite{lauter:2009}, a cryptographic hash function was defined: $$h:\{ 0,1\}^r \rightarrow \{ 0,1\}^s$$ based on the Supersingular Isogeny Graph (SSIG) for a fixed prime $p$ of cryptographic size, and a fixed small prime $\ell \ne p$.  The hash function processes the input string in blocks which are used as directions for walking around the graph starting from a given fixed vertex.  The output of the hash function is the $j$-invariant of an elliptic curve over $\FF_{p^2}$ which requires $2\log(p)$ bits to represent, so $m= 2 \lceil \log(p)\rceil$.
For the security of the hash function, it is necessary to avoid the generic {\it birthday attack}. This attack runs in time proportional to the square root of the size of the graph, which is the {\it Eichler class number}, roughly $\lfloor p/12 \rfloor$. So in practice, we must pick $p$ so that $\log(p) \approx 256$.

The integer $r$ is the length of the bit string input to the hash function. If $\ell = 2$, which is the easiest case to implement and a common choice, then $r$ is precisely the number of steps taken on the walk in the graph, since the  graph is 3-regular, with no backtracking allowed, so the input is processed bit-by-bit.   In order to assure that the walk reaches a sufficiently random vertex in the graph, the number of steps should be roughly $\log(p) \approx 256$.  A CGL-hash function is thus specified by giving the primes $p$, $\ell$, the starting point of the walk, and the integers $r \approx 256$, $s$.  (Extra congruence conditions were imposed on $p$ to make it an undirected graph with no small cycles.)

The hard problems stated in~\cite{lauter:2009} corresponded to the important security properties of {\it collision} and {\it preimage resistance} for this hash function.  For preimage resistance, the problem~\cite[Problem 3]{lauter:2009} stated was: given $p$, $\ell$, $r>0$, and two supersingular $j$-invariants modulo $p$, to find a path of length $r$ between them:

\begin{problem}(Path-finding \cite{lauter:2009})
\label{pbm:path}
 Let $p$ and $\ell$ be distinct prime numbers, $r> 0$, and $E_0$ and $E_1$ two supersingular elliptic curves over $\mathbb{F}_{p^2}$.
 Find a path of length $r$ in the $\ell$-isogeny graph corresponding to a composition of $r$ $\ell$-isogenies leading from $E_0$ to $E_1$ (i.e. an isogeny of degree $\ell^r$ from $E_0$ to $E_1$).
\end{problem}

It is worth noting that, to break the preimage resistance of the specified hash function, you must find a path of exactly length $r$, and this is analogous to the situation for breaking the security of the key-exchange protocol.  However, the problem of finding *any* path between two given vertices in the SSIG graphs is also still open.  For the LPS graphs, the algorithm presented in~\cite{SCN:PetLauQui08} did not find a path of a specific given length, but it was still considered to be a ``break'' of the hash function.

Furthermore, the diameter of these graphs, both LPS and SSIG graphs, has been extensively studied. It is known that the diameter of the graphs is roughly $\log(p)$ (it is $c\log(p)$, where $c$ is a constant between $1$ and $2$, (see for example~\cite{S17})).  That means that if $r$ is greater than $c\log(p)$, then given two vertices, it is likely that a path of length $r$ between them may  exist.  The fact that walks of length greater than $c\log(p)$ approximate the uniform distribution very closely means that you are not likely to miss any significant fraction of the vertices with paths of that length, because that would constitute a bias. Also, if $r \gg \log(p)$ then there may be many paths of length $r$.  However, if $r$ is much less than $\log(p)$, such as $\frac{1}{2}\log(p)$, there may be {\it no path} of such a short length between two given vertices.  See~\cite{LP15} for a discussion of the ``sharp cutoff'' property of Ramanujan graphs.

But in the cryptographic applications, given an instance of the key-exchange protocol to be attacked, we {\it know} that there exists a path of length $n$ between $E$ and $E_A$, and the hard problem is to find it.  
The set-up for the key-exchange requires $p=\ell_A^n \ell_B^m \pm 1$, where $n$ and $m$ are roughly the same size, and $\ell_A$ and $\ell_B$ are very small, such as $\ell_A = 2$ and $\ell_B = 3$.  It follows that $n$ and $m$ are both approximately half the diameter of the graph (which is roughly $\log(p)$).  So it is unlikely to find paths of length $n$ or $m$ between two random vertices.  If a path of length $n$ exists and Algorithm A finds a path, then it is very likely to be the one which was constructed in the key exchange.  If not, then Algorithm A can be repeated any constant number of times.  So we have the following reduction:

\begin{theorem}
\label{lemma::reduction}
Assume as for the Key Exchange set-up that $p =\ell_A^n\cdot \ell_B^m + 1$ is a prime of cryptographic size, i.e. $\log(p) \ge 256$, $\ell_A$ and $\ell_B$ are small primes, such as $\ell_A = 2$ and $\ell_B = 3$, and $n \approx m$ are approximately equal.
Given an algorithm to solve Problem \ref{pbm:path} (Path-finding), it can be used to solve Problem \ref{pbm:sscdh} (Key Exchange) with overwhelming probability.  The failure probability is roughly $$\frac{\ell_A^n + \ell_A^{n-1}}{p} \approx \frac{\sqrt{p}}{p}.$$
\end{theorem}

 \begin{proof}
Given an algorithm (Algorithm A) to solve Problem~\ref{pbm:path}, we can use this to solve Problem~\ref{pbm:sscdh} as follows. Given $E$ and $E_A$, use Algorithm A to find the path of length $n$ between these two vertices in the $\ell_A$-isogeny graph.   Now use Lemma \ref{lemma::isogeny-composition-2} below to produce a point $R_A$ which generates the $\ell_A^n$-isogeny between $E$ and $E_A$. Repeat this to produce the point $R_B$ which generates the $\ell_B^m$-isogeny between $E$ and $E_B$ in the $\ell_B$-isogeny graph.  Because the subgroups generated by $R_A$ and $R_B$ have smooth order, it is easy to write $R_A$ in the form $[m_A]P_A + [n_A]Q_A$ and $R_B$ in the form $[m_B]P_B + [n_B]Q_B$.
 Using the knowledge of $m_A$, $n_A$, $m_B$, $n_B$, we can construct $E_{AB}$ and recover the $j$-invariant of $E_{AB}$, allowing us to solve Problem~\ref{pbm:sscdh}.

The reason for the qualification ``with overwhelming probability'' in the statement of the theorem is that it is possible that there are multiple paths of the same length between two vertices in the graph. If there are multiple paths of length $n$ (or $m$) between the two vertices, it suffices to repeat Algorithm A to find another path.  This approach is sufficient to break the Key Exchange if there are only a small number of paths to try.  As explained above, with overwhelming probability, there are {\it no} other paths of length $n$ (or $m$) in the Key Exchange setting.  

In the SSIG corresponding to $(p, \ell_A)$, the vertices $E$ and $E_A$ are a distance of $n$ apart.  Starting from the vertex $E$ and considering all paths of length $n$, the number of possible endpoints is at most $\ell_A^n + \ell_A^{n-1}$ (See Corollary~\ref{corollary::number-isogenies} below).  Considering that the number of vertices in the graph is roughly $\lfloor p/12 \rfloor$, then the probability that a given vertex such as $E_A$ will be the endpoint of one of the walks of length $n$ is roughly
$$\frac{\ell_A^n + \ell_A^{n-1}}{p} \approx \frac{\sqrt{p}}{p} \le 2^{-128}.$$
This estimate does not use the Ramanujan property of the SSIG graphs.  While a generic random graph could potentially have a topology which creates a bias towards some subset of the nodes, Ramanujan graphs cannot, as shown in~\cite[Theorem 3.5]{LP15}.
 \end{proof}

\section{Composing isogenies}\label{sec::correspondence}

Let $k$ be a positive integer. Every separable $k$-isogeny $\phi : E_0 \rightarrow E_1$ is determined by its kernel up to composition with an automorphism of the elliptic curve $E_1.$ Thus the edge corresponding to $\phi$ is uniquely determined by $\ker(\phi )$ and vice versa. This kernel is a subgroup of the $k$-torsion $E_0[k]$, and the latter is isomorphic to $\ZZ/k\ZZ \times \ZZ/k\ZZ$ if $k$ is coprime to the characteristic of the field we are working over.

Hence, fixing a prime $\ell$ and working over a finite field $\FF_q$ which has characteristic different from $\ell$, the number of $\ell$-isogenies $\phi : E_0 \rightarrow E_1$ that correspond to different edges of the graph is equal to the number of subgroups of $\ZZ/\ell\ZZ \times \ZZ/\ell\ZZ$ of order $\ell$. It is well known that this number is equal to $\ell+1$. In other words, $E$ is $\ell$-isogenous to precisely $\ell+1$ elliptic curves.

However, some of these $\ell$-isogenous curves may be isomorphic. Therefore, in the isogeny graph (where nodes represent isomorphism classes of curves), $E$ has degree $\ell+1$ and may have $\ell+1$ neighbors or fewer.

Using V\'{e}lu's formulas, the equations for an edge can be computed from its kernel. Hence for computational purposes, it is important to write down this kernel explicitly. This is best done by specifying generators. Let $P, Q \in E_0$ be the generators of $E_0[\ell] \cong \ZZ/\ell\ZZ \times \ZZ/\ell\ZZ$. Then the subgroups of order $\ell$ are generated by $Q$ and $P + iQ$ for $i = 0,\ldots,\ell-1$.

We now study isogenies obtained by composition, and isogenies of degree a prime power. It turns out that these correspond to each other under certain conditions. The first condition is that the isogeny is cyclic. Notice that every prime order group is cyclic, therefore all $\ell$-isogenies are cyclic (meaning they have cyclic kernel). However, this is not necessarily true for isogenies whose order is not a prime. The second condition is that there is no backtracking, defined as follows:

\begin{definition}
 For a chain of isogenies $\phi_m \circ \phi_{m-1} \circ \ldots \circ \phi_1$ ($\phi_i:E_{i-1}\rightarrow E_i$), we say that it has \emph{no backtracking} if $\phi_{i+1}\ne \epsilon\circ \hat\phi_i $ for all $i = 1,\ldots,m-1$ and any $\epsilon\in \Aut (E_{i+1}) $, since this corresponds to a walk in the $\ell$-isogeny graph without backtracking.
\end{definition}

In the following, we show that chains of $\ell$-isogenies of length $m$ without backtracking correspond to cyclic $\ell^m$-isogenies. Recall that we are only considering separable isogenies throughout.

\begin{lemma}\label{lemma::isogeny-composition-1}
Let $\ell$ be a prime, and let $\phi$ be a separable $\ell^m$-isogeny with cyclic kernel. Then there exist cyclic $\ell$-isogenies $\phi_1,\ldots,\phi_{m}$ such that
$ \phi = \phi_m \circ \phi_{m-1} \circ \ldots \circ \phi_1 $ without backtracking.
\end{lemma}

\begin{proof}
 Assume that $\phi = E_0 \rightarrow E$, and that its kernel is $\langle P_0 \rangle \subseteq E_0$, where $P_0$ has order $\ell^m$. For $i = 1,\ldots,m$, let
 $$ \phi_i : E_{i-1} \rightarrow E_i $$
 be an isogeny with kernel $\langle \ell^{m-i}P_{i-1} \rangle$, where $P_i = \phi_{i}(P_{i-1}).$

 We show that $\phi_i$ is an $\ell$-isogeny for $i \in \{1,\ldots,m\}$ by observing that $\ell^{m-i}P_{i-1}$ has order $\ell$. The statement is trivial for $i=1$. For $i \geq 2$, clearly $\ell^{m-i}P_{i-1} = \ell^{m-i}\phi_{i-1}(P_{i-2}) = \phi_{i-1}(\ell^{m-i}P_{i-2}) \neq \O$, since $\ell^{m-i}P_{i-2} \notin \ker \phi_{i-1} = \langle \ell^{m-(i-1)}P_{i-2} \rangle = \{ \ell^{m-(i-1)}P_{i-2}, 2\ell^{m-(i-1)}P_{i-2},\ldots,(\ell-1)\ell^{m-(i-1)}P_{i-2} \} $.  Furthermore, $\ell \cdot \ell^{m-i}P_{i-1} = \ell^{m-(i-1)}\phi_{i-1}(P_{i-2}) = \phi_{i-1}(\ell^{m-(i-1)}P_{i-2}) = \O$, using the definition of $\ker \phi_{i-1}$.

 Next, we show by induction that $\phi_i \circ \ldots \circ \phi_1$ has kernel $\langle \ell^{m-i}P_0 \rangle$. Then it follows that $\phi_m \circ \ldots \circ \phi_1$ is the same as $\phi$ up to an automorphism $\epsilon$ of $E$, since the two have the same kernel. Replacing $\phi_m$ with $\epsilon\circ \phi_m$ if necessary we have $ \phi = \phi_m \circ \phi_{m-1} \circ \ldots \circ \phi_1 .$
 The case $i = 1$ is trivial: $\phi_1 : E_0 \rightarrow E_1$ has kernel $\langle \ell^{m-1}P_0 \rangle$ by definition. Now assume the statement is true for $i-1$. Then, we have $\langle \ell^{m-i}P_0 \rangle \subseteq \ker \phi_{i} \circ \ldots \circ \phi_1$. Conversely, let $Q \in \ker \phi_i\circ\ldots\circ\phi_1$. Then $\phi_{i-1}\circ\ldots\circ\phi_i(Q) \in \ker\phi_i = \langle \ell^{m-i}P_{i-1} \rangle = \phi_{i-1}(\langle \ell^{m-i}P_{i-2} \rangle) = \ldots = \phi_{i-1}\circ\ldots\circ\phi_1(\langle \ell^{m-i}P_0 \rangle)$ and hence $Q \in \langle \ell^{m-i}P_0 \rangle + \ker\phi_{i-1}\circ\ldots\circ\phi_1 = \langle \ell^{m-i}P_0 \rangle + \langle \ell^{m-(i-1)}P_0 \rangle = \langle \ell^{m-i}P_0 \rangle$.

 Finally, we show that there is no backtracking in $\phi_m\circ\ldots\circ\phi_1$. Contrarily, assume that there is an $i \in \{1,\ldots,m-1\}$ and $\epsilon\in \Aut(E_{i+1})$ such that $\phi_{i+1}=\epsilon\circ \hat \phi_i$. Then, since $\ker (\phi_{i+1}\circ \phi_i )=\ker(\epsilon\circ \hat\phi_i\circ \phi_i )= \ker [\ell]$, we have $\ker(\phi_{i+1}\circ \phi_i \circ\phi_{i-1}\circ\ldots\circ\phi_1) =\ker ([\ell]\circ\phi_{i-1}\circ\ldots\circ\phi_1)$. Notice that $[\ell]$ commutes with all $\phi_j$, and hence $E_0[\ell] \subseteq \ker(\phi_{i+1}\circ \phi_i \circ\phi_{i-1}\circ\ldots\circ\phi_1)\subseteq \ker(\phi_{m}\circ \phi_i \circ\phi_{i-1}\circ\ldots\circ\phi_1)=\ker\phi$. Since $E_0[\ell] \cong \ZZ/\ell\ZZ \times \ZZ/\ell\ZZ$, the kernel of $\phi$ cannot be cyclic, a contradiction.
\end{proof}

\begin{remark}
  It is clear that in the above lemma, if $\phi$ is defined over a finite field $\FF_q$, then all $\phi_i$ are also defined over this field. Namely, if $E_0$ is defined over $\FF_q$ and the kernel is generated by an $\FF_q$-rational point, then by V\'elu we obtain $\FF_q$-rational formulas for $\phi_1$, which means that $\phi_1$ is defined over $\FF_q$, and so on.
\end{remark}

\begin{lemma}\label{lemma::isogeny-composition-2}
 Let $\ell$ be a prime, let $E_i$ be elliptic curves for $i = 0,\ldots,m$, and let $\phi_i : E_{i-1} \rightarrow E_i$ be $\ell$-isogenies for $i = 1,\ldots,m$ such that $\phi_{i+1}\ne \epsilon\circ \hat\phi_i $ for $i = 1,\ldots,m-1$ and any $\epsilon\in \Aut (E_{i+1}) $ (i.e.\ there is no backtracking). Then $\phi_m \circ \ldots \circ \phi_1$ is a cyclic $\ell^m$-isogeny.
\end{lemma}

\begin{proof}
 The degree of isogenies multiplies when they are composed, see e.g.\ \cite[Ch.\ III.4]{silverman}. Hence we are left with proving that the composition of the isogenies is cyclic.

 First note that all $\phi_i$ are cyclic since they have prime degree, and denote by $P_{i-1} \in E_{i-1}$ the generators of the respective kernels. Let $Q_{m-1}$ be a point on $E_{m-1}$ such that $\ell Q_{m-1} = P_{m-1}$. Notice that such a point always exists over the algebraic closure of the field of definition of the curve. Let $R_{m-2} = \hat{\phi}_{m-1}(Q_{m-1})$, where the hat denotes the dual isogeny. Then $\phi_m \circ \phi_{m-1}(R_{m-2}) = \phi_m \circ \phi_{m-1} \circ \hat{\phi}_{m-1}(Q_{m-1}) = \phi_m \circ [\ell] (Q_{m-1}) = \phi_m(\ell Q_{m-1}) = \phi_m(P_{m-1}) = \mathcal{O}$, and hence $R_{m-2}$ is in the kernel of $\phi_m \circ \phi_{m-1}$.

 Next we show that $R_{m-2}$ has order $\ell^2$, which implies that it generates the kernel of $\phi_m \circ \phi_{m-1}$. Suppose that $\ell R_{m-2} = \O$. Then $\O = \ell R_{m-2} = \ell \hat{\phi}_{m-1}(Q_{m-1}) = \hat{\phi}_{m-1}(P_{m-1})$. Since $P_{m-1}$ has order $\ell$, this implies that $P_{m-1}$ generates the kernel of $\hat{\phi}_{m-1}$. However, $P_{m-1}$ also generates the kernel of $\phi_m$, so $\epsilon\circ \hat{\phi}_{m-1} = \phi_m$ for some $\epsilon\in \Aut(E_m)$. But this is a contradiction to the assumption of no backtracking.

 By iterating this argument, we obtain a point $R_0$ which generates the kernel of $\phi_m\circ\ldots\circ\phi_1$, and hence this isogeny is cyclic.
\end{proof}

Combining Lemmas \ref{lemma::isogeny-composition-1} and \ref{lemma::isogeny-composition-2}, we obtain the following correspondence.

\begin{corollary}\label{corollary::correspondence}
 Let $\ell$ be a prime and $m$ a  positive integer. There is a one-to-one correspondence between cyclic separable $\ell^m$-isogenies and chains of separable $\ell$-isogenies of length $m$ without backtracking. (Here we do not distinguish between isogenies that differ by composition with an automorphism on the image.)
\end{corollary}

Next, we investigate how many such isogenies there are. We start by studying $\ell^m$-isogenies. The following group theory result is crucial.

\begin{lemma}\label{lemma::any-subgroups}
 Let $\ell$ be a prime and $m$ a positive integer. Then the number of subgroups of $\ZZ/\ell^m\ZZ \times \ZZ/\ell^m\ZZ$ of order $\ell^m$ is $\frac{\ell^{m+1}-1}{\ell-1}$, and $\ell^m + \ell^{m-1}$ of these subgroups are cyclic.
\end{lemma}

\begin{proof}
Every subgroup of $\ZZ/\ell^m\ZZ \times \ZZ/\ell^m\ZZ$
is isomorphic to $\ZZ/\ell^i\ZZ \times \ZZ/\ell^j\ZZ$ for $0 \leq i \leq j \leq m$. The number of subgroups
which are isomorphic to $\ZZ/\ell^i\ZZ \times \ZZ/\ell^j\ZZ$ is 1 if $i=j$ and $\ell^{j-i} + \ell^{j-i-1}$
otherwise.

A direct consequence of the above statement is that there are
\[  \sum_{i=0}^{\floor{\frac{m-1}{2}}} \ell^{m-2i} + \ell^{m-2i-1} + \epsilon_m = \sum_{t=0}^{m} \ell^t    \]
subgroups, where $\epsilon_m = 0$ if $k$ is odd and 1 otherwise. This proves the first statement.

For the second statement, let $H$ be a cyclic subgroup of $\ZZ/\ell^m\ZZ \times \ZZ/\ell^m\ZZ$
of order $l^m$. Then $H$ is generated by an element of $\ZZ/\ell^m\ZZ \times \ZZ/\ell^m\ZZ$
of order $l^m$, and contains $l^m - l^{m-1}$ elements of order $l^m$. Therefore, the number
of such subgroups is the number of elements of $\ZZ/\ell^m\ZZ \times \ZZ/\ell^m\ZZ$ of order
$l^m$ divided by $l^m - l^{m-1}$.

Let $(a,b)$ be an element of $\ZZ/\ell^m\ZZ \times \ZZ/\ell^m\ZZ$ of order $l^m$. Then one of
$a$ or $b$ has order $l^m$. If $a$ has order $l^m$, then there are $\varphi(\ell^m)=l^m - l^{m-1}$ choices for
$a$, and $l^m$ for $b$. That is, there are $l^m\cdot(l^m-l^{m-1})$ choices in total.

Otherwise, there are $l^{m-1}$ choices for $a$ (representing the number of elements of
order at most $l^{m-1}$), and $l^m - l^{m-1}$ choices for $b$. That is, there are
$l^{m-1}\cdot(l^m-l^{m-1})$ choices in total. This means the total number of cyclic
subgroups of  $\ZZ/\ell^m\ZZ \times \ZZ/\ell^m\ZZ$ of order $l^m$ is

\[ \frac{l^m\cdot(l^m-l^{m-1}) + l^{m-1}\cdot(l^m-l^{m-1})}{l^m - l^{m-1}}= l^m + l^{m-1}. \]
\end{proof}

\begin{remark}
One could also see the first statement in the lemma above by noting that this is the same as the degree of the Hecke operator $T_{\ell^m}$ which is $\sigma_1(\ell^m)$. We thank the referee for pointing this out.
\end{remark}

\begin{corollary}\label{corollary::number-isogenies}
There are $\frac{\ell^{m+1}-1}{\ell-1}$ separable $\ell^m$-isogenies originating at a fixed elliptic curve, and $\ell^m + \ell^{m-1}$ of them are cyclic. (Here we are counting isogenies as different if they differ even after composition with any automorphism of the image.)
\end{corollary}

Using the correspondence from Corollary \ref{corollary::correspondence}, we then obtain the following.

\begin{theorem}\label{thm::number-chains}
 The number of chains of $\ell$-isogenies of length $m$ without backtracking is $\ell^m+ \ell^{m-1}$. (Here we do not distinguish between isogenies that differ by composition with an automorphism on the image.)
\end{theorem}

This last result can be observed in a much more elementary way, which is also enlightening. We consider chains of $\ell$-isogenies of length $m$. To analyze the situation, it is helpful to draw a graph similar to an $\ell$-isogeny graph but that does {\it not} identify isomorphic curves.
This graph is an $(\ell+1)$-regular tree of depth $m$. The root of the tree has $\ell+1$ children, and every other node (except the leaves) has $\ell$ children. The leaves have depth $m$. It is easy to work out that the number of leaves in this tree is $(\ell+1)\ell^{m-1}$, and this is also equal to the number of paths of length $m$ without backtracking, as stated in Theorem \ref{thm::number-chains}.

Finally, this graph also helps us count the number of chains of $\ell$-isogenies of length $m$ including those that backtrack. By examining the graph carefully, we can see that the number of such walks is $\ell^m + \ell^{m-1} + \ldots + \ell + 1$, and according to Corollary \ref{corollary::number-isogenies}, this corresponds to the number of $\ell^m$-isogenies that are not necessarily cyclic.

These results were also observed experimentally using Sage. The numbers match the results of our experiments for small values of $\ell$ and $m$, over various finite fields and for different choices of elliptic curves, see Table \ref{table::endpoints}. Notice that the images under isogenies with distinct kernels may be isomorphic, leading to double edges in an isogeny graph that identifies isomorphic curves. Hence, the number of isomorphism classes of images (i.e.\ the number of neighbors in the isogeny graph) may be smaller than the number of isogenies stated in the table. 

\begin{table}[h]
\begin{center}
\begin{tabular}{| c | c | c | c |}
 \hline
 $\ell$ & $m$ & number of isogenies & number of isogenies  \\
  &  & without backtracking & with backtracking \\\hline
 2 & 4 & 24 & 31 \\
 2 & 5 & 48 & 63 \\
 2 & 6 & 96 & 127 \\
 2 & 7 & 192 & 255 \\
 3 & 4 & 108 & 121 \\
 3 & 5 & 324 & 364 \\\hline
\end{tabular}
\end{center}
\caption{For small fixed $\ell$ and $m$, values obtained experimentally for the number of $\ell$-isogeny-chains of length $m$ starting at a fixed elliptic curve $E$ without and with backtracking.}
\label{table::endpoints}
\end{table}

\part{\Large Constructions of Ramanujan graphs}\label{part2}

\newcommand{\ve}{{\mathbf{e}}}
\newcommand{\vv}{{\mathbf{v}}}
\newcommand{\vu}{{\mathbf{u}}}
\newcommand{\GLt}{{\mathrm{GL}}_2}
\newcommand{\SLt}{{\mathrm{SL}}_2}

\newcommand{\Qv}{{\QQ_v}}
\newcommand{\Zv}{{\ZZ_v}}
\newcommand{\Qp}{{\QQ_p}}
\newcommand{\Zp}{{\ZZ_p}}

In this section we review the constructions of two families of Ramanujan graph, LPS graphs and Pizer graphs. Ramanujan graphs are optimal expanders; see Section \ref{subset:autoprelim} for some related background. The purpose is twofold. On the one hand we wish to explain how equivalent constructions on the same object highlight different significant properties. On the other hand, we wish to explicate the relationship between LPS graphs and Pizer graphs.

Both families (LPS and Pizer) of Ramanujan graphs can be viewed (cf.\ \cite[Section 3]{li:1996}) as a set of ``local double cosets'', i.e.\ as a graph on
\begin{equation}\label{eq:doublecoset-PGL}
\Gamma \backslash \PGLt(\Ql)/\PGLt(\Zl),
\end{equation}
where $\Gamma $ is a discrete cocompact subgroup.
In both cases, one has a chain of isomorphisms that are used to show these graphs are Ramanujan, and in both cases one may in fact vary parameters to get an infinite family of Ramanujan graphs.

To explain this better, we introduce some notation. Let us choose a pair of distinct primes $p$ and $l$ for an $(l+1)$-regular graph whose size depends on $p.$ (An infinite family of Ramanujan graphs is formed by varying $p.$) Let us fix a quaternion algebra $B$ defined over $\QQ$ and ramified at exactly one finite prime and at $\infty ,$ and an order of the quaternion algebra $\O .$ Let $\A$ denote the ad\`{e}les of $\QQ$ and $\A_f$ denote the finite ad\`{e}les. For precise definitions see Section \ref{subset:autoprelim}.

In the case of Pizer graphs, let $B=B_{p,\infty}$  be ramified at $p$ and $\infty ,$ and take $\O$ to be a maximal order (i.e.\ an order of level $p$).\footnote{A similar construction exists for a more general $\O .$ However, to relate the resulting graph to supersingular isogeny graphs, we require $\O$ to be maximal. } Then we may construct (as in \cite{pizer:1998}) a graph by giving its adjacency matrix as a Brandt matrix. (The Brandt matrix is given via an explicit matrix representation of a Hecke operator associated to $\O$.) Then we have (cf.\ \cite[(1)]{lauter:2009r}) a chain of isomorphisms connecting \eqref{eq:doublecoset-PGL} with supersingular isogeny graphs (SSIG) discussed in Part \ref{part1} above:
\begin{equation}\label{eq:whatEyalwrote}
(\O[l^{-1}])^{\times }\backslash \GLt(\Ql)/\GLt(\Zl) \cong B^{\times }(\QQ)\backslash B^{\times }(\A_f)/B^{\times }(\hat{\ZZ})\cong \Cl \O \cong \text{SSIG}.
\end{equation}

This can be used (cf.\ \cite[5.3.1]{lauter:2009r}) to show that the supersingular $l$-isogeny graph is connected, as well as the fact that it is indeed a Ramanujan graph.

In the case of LPS graphs the choices are very different. Let $B=B_{2,\infty }$ now be the Hamiltonian quaternion algebra. The group $\Gamma $ in \eqref{eq:doublecoset-PGL} is chosen as a congruence subgroup dependent on $p.$ This leads to a larger graph whose constructions fits into the following chain of isomorphisms:
\begin{equation}\label{eq:LPSchain1}
\PSLt(\FF_p)\cong \Gamma(2p)\backslash \Gamma(2)\cong \Gamma(2p)\backslash T\cong \Gamma(2p)\backslash \PGLt(\Ql)/\PGLt(\Zl)\cong G'(\QQ)\backslash H_{2p}/G'(\RR)K_0^{2p}.
\end{equation}
The isomorphic constructions and their relationship will be made explicit in Sections \ref{subsubsect:CaylePres}-\ref{subsubsect:localdoublecosets} and Section \ref{subsubsect:SA_LPS}.
We shall also explain how properties of the graph, such as its regularity, connectedness and the Ramanujan property, are highlighted by this chain of isomorphisms. For now we give only an overview, to be able to compare this case with that of Pizer graphs. The quotient $\PGLt(\Ql)/\PGLt(\Zl)$ has a natural structure of an infinite tree $T.$ This tree can be defined in terms of homothety classes of rank two lattices of $\Ql\times \Ql$ (see Section \ref{subsubsect:infinitelatticetree}). One may define a group $G'=B^{\times }/Z(B^{\times })$ and its congruence subgroups $\Gamma (2)$ and $\Gamma (2p),$ and show that the discrete group $\Gamma (2)$ acts simply transitively on the tree $T,$ and hence $\Gamma(2p)\backslash T$ is isomorphic to the finite group $ \Gamma(2)/\Gamma(2p) .$ Using the Strong Approximation theorem, this turns out to be isomorphic to the group $\PSLt(\FF_p) .$ The latter has a structure of an $(l+1)$-regular Cayley graph. A second application of the Strong Approximation Theorem with $K_0^{2p}$, an open compact subgroup of $G'(\A_f)$, shows that  $H_{2p}$ is a finite index normal subgroup of $G'(\A)$.

Note that an immediate distinction between Pizer and LPS graphs is that the quaternion algebras underlying the constructions are different: they ramify at different finite primes ($p$ and $2,$ respectively). In addition, the size of the discrete subgroup $\Gamma $ determining the double cosets of \eqref{eq:doublecoset-PGL} is different in the two cases. Accordingly, the size of the resulting graphs is different as well. We shall see
that (under appropriate assumptions on $p$ and $l$) the Pizer graph has $\frac{p-1}{12}$ vertices, while the LPS graph has order $|\PSLt(\FF_p)|=\frac{p(p^2-1)}{2}.$ One may consider an order $\O_{LPS}$ such that $(\O_{LPS}[l^{-1}])^{\times }\cong \Gamma (2p)$
analogously to the relationship of $\O $ and $\Gamma $ in the Pizer case and \eqref{eq:whatEyalwrote}. However, this order $\O_{LPS}$ is unlike the Eichler order from the Pizer case. (It has a much higher level.) In particular, there is a discrepancy between the order of the class set $\Cl \O_{LPS}$ and the order of the LPS graph. This is a numerical obstruction indicating that an analogue of the chain \eqref{eq:whatEyalwrote} for LPS graphs is at the very least not straightforward.

The rest of the paper has the following outline. In Section \ref{subsect:LPS} we explore the isomorphic constructions of LPS graphs from \eqref{eq:LPSchain1}. We give the construction as a Cayley graph in Section \ref{subsubsect:CaylePres}. The infinite tree of homothety classes of lattices is given in Section \ref{subsubsect:infinitelatticetree}. In Section \ref{subsubsect:localdoublecosets} we explain how local double cosets of the Hamiltonian quaternion algebra connect these constructions. Section \ref{app:sect:explicitisom} makes one step of the chain of isomorphisms in \eqref{eq:LPSchain1} completely explicit in the case of $l=5$ and $l=13,$ and describes how the same can be done in general. In Section \ref{subsect:StrongApprox} we give an overview of how Strong Approximation plays a role in proving the isomorphisms and the connectedness and Ramanujan property of the graphs. In Section \ref{sect:Pizer} we turn briefly to Pizer graphs. We summarize the construction, and explain how various restrictions on the prime $p$ guarantee properties of the graph. Section \ref{subsect:Pizer_primes} contains the computation of a prime $p$ where the existence of both an LPS and a Pizer construction is guaranteed (for $l=5$). In Section \ref{subsect:LPS_Pizer_compare} we say a bit more of the relationship of Pizer and LPS graphs, having introduced more of the objects mentioned in passing above.

Throughout this part of the paper we aim to only include technical details if we can make them fairly self-contained and explicit, and otherwise to give a reference for further information.

\section{Background on Ramanujan graphs  and ad\`{e}les}\label{subset:autoprelim}
In this section we fix notation and review some definitions and facts that we will be using for the remainder of Part \ref{part2}.

Expander graphs are graphs where small sets of
vertices have many neighbors. For many applications of expander graphs, such as in Part \ref{part1}, one wants $(l+1)$-regular expander graphs $X$ with $l$ small and the number of vertices of $X$ large.
If $X$ is an $(l+1)$-regular graph (i.e. where every vertex has degree $l+1$), then $l+1$ is an eigenvalue of the adjacency matrix of $X.$ All eigenvalues $\lambda $ satisfy $-(l+1)\leq \lambda \leq (l+1)$, and $-(l+1)$ is an eigenvalue if and only if $X$ is bipartite.
Let $\lambda(X)$ be the second largest eigenvalue in absolute value of the adjacency matrix. The smaller $\lambda(X)$ is, the better expander $X$ is. Alon--Boppana proved that for  an \emph{infinite} family of $(l+1)$-regular graphs of increasing size,  $\lim \inf_{(X)}\lambda (X)\geq 2 \sqrt{l}$  \cite{alon:1986}. An $(l+1)$-regular graph $X$ is called Ramanujan if $\lambda(X)\leq 2 \sqrt l$.  Thus an infinite family of Ramanujan graphs are optimal expanders.

For a finite prime $p$,  let $\QQ_p$ denote the field of $p$-adic numbers and $\ZZ_p$ its ring of integers. Let $\QQ_\infty=\RR$. We denote the ad\`{e}le ring of $\QQ$ by $\A$ and recall that it is defined as a restricted direct product in the following way,
\[
\A=\sideset{}{'}\prod_p\QQ_p=\left\{ (a_p)\in \prod_p \QQ_p : a_p \in \ZZ_p \text{ for  all but a finite number of  $p< \infty$}\right\}.
\]
We denote the ring of finite ad\`{e}les by $\A_f$, that is
\[
\A_f=\sideset{}{'}\prod_{p<\infty}\QQ_p=\left\{ (a_p)\in \prod_{p<\infty} \QQ_p : a_p \in \ZZ_p \text{ for  all but a finite number of  $p$}\right\}.
\]

Let $\A^{\times}$ denote the id\`{e}le group of $\QQ$, the group of units of $\A$,
\[
\A^{\times}=\sideset{}{'}\prod_p\QQ_p=\left\{ (a_p)\in \prod_p \QQ_p^\times : a_p \in \ZZ_p^\times \text{ for  all but a finite number of $p< \infty$}\right\}.
\]

Let $B$ be a quaternion algebra over $\QQ$, $B^\times$ the invertible elements of $B$ and $\O$ an order of $B$. For a prime $p$ let $\O_p=\O\otimes_{\ZZ}\ZZ_p$. Then let
\[
B^\times(\A)=\sideset{}{'}\prod_pB^\times(\QQ_p)=\left\{(g_p)\in \prod_p B^\times(\QQ_p) : g_p\in \O_p^\times \text{ for  all  but a finite number of $p< \infty$}\right\}.
\]

More generally for an indexed set of locally compact groups $\{G_v\}_{v\in I}$ with a corresponding indexed set of compact open subgroups $\{K_v\}_{v\in I}$ we may define the restricted direct product of the $G_v$ with respect to the $K_v$ by the following
\[G:=\sideset{}{'}\prod_{v\in I} G_v=\left\{ (g_v)\in \prod_{v\in I} G_v : g_v \in K_v \text{ for  all  but a finite number of $v$} \right\}.
\]
If we define  a neighborhood base of the identity as
\[
\left\{\prod_v U_v : U_v \text{ neighborhood of identity in $G_v$ and $U_v=K_v$ for  all but a finite number of $v$}\right\}
\]
then $G$ is a locally compact topological group.

\section{LPS Graphs}\label{subsect:LPS}

We describe the LPS graphs used in \cite{lauter:2009} for a proposed hash function. They were first considered in \cite{lubotzky:1988}, for further details see also \cite{lubotzky:2010}. We shall examine the objects and isomorphisms in \eqref{eq:LPSchain1} in more detail. We review constructions of these graphs in turn as Cayley graphs and graphs determined by rank two lattices or, equivalently, local double cosets. Throughout this section, let $l$ and $p$ be distinct, odd primes both congruent to $1$ modulo $4.$ We shall give constructions of $(l+1)$-regular Ramanujan graphs whose size depends on $p.$
We shall also assume for convenience\footnote{If $p$ is not a square modulo $l,$ then the constructions described below result in bipartite Ramanujan graphs with twice as many vertices.} that $\left(\frac{p}{l}\right)=1,$ i.e. that $p$ is a square modulo $l.$

\subsection{Cayley graph over $\FF_p.$}\label{subsubsect:CaylePres}
This description follows \cite[Section 2]{lubotzky:1988}. The graph we are interested in is the Cayley graph of the group $\PSLt(\FF_p).$ We specify a set of generators $S$ below. The vertices of the graph are the $\frac{p(p^2-1)}{2}$ elements of $\PSLt(\FF_p).$ Two vertices $g_1,g_2\in \PSLt(\FF_p)$ are connected by an edge if and only if $g_2=g_1h$ for some $h\in S.$

Next we give the set of generators $S.$ Since $l\equiv 1\mod{4}$ it follows from a theorem of Jacobi \cite[Theorem 2.1.8]{lubotzky:2010} that there are $l+1$ integer solutions to
\begin{equation}\label{eq:lsquaresumconds}
l=x_0^2+x_1^2+x_2^2+x_3^2;\ \ 2\nmid x_0;\ \ x_0>0.
\end{equation}
In this case we will also have $2|x_i$ for all $i>0.$ 
Let $S$ be the set of solutions of \eqref{eq:lsquaresumconds}. Since $p\equiv1\mod{4}$ we have $\left(\frac{-1}{p}\right)=1.$ Let $\varepsilon \in \ZZ $ such that $\varepsilon^2\equiv -1\mod{p}.$ Then to each solution of \eqref{eq:lsquaresumconds} we assign an element of $\PGLt(\ZZ)$ as follows:
\begin{equation}\label{eq:elementsofSmappedtomtcesmodp}
(x_0,x_1,x_2,x_3)\mapsto \tbtmtx{x_0+x_1\varepsilon}{x_2+x_3\varepsilon}{-x_2+x_3\varepsilon}{x_0-x_1\varepsilon}.
\end{equation}
Note that the matrix on the right-hand side has determinant $l\mod{p}.$ Since $\left(\frac{l}{p}\right)=1$ this determines an element of $\PSLt(\FF_p).$
The $l+1$ elements of $\PSLt(\FF_p)$ determined by $\eqref{eq:elementsofSmappedtomtcesmodp}$ form the set of Cayley generators. Let us abuse notation and denote this set with $S$ as well.
This graph is connected. To prove this fact, one may use the theory of quadratic Diophantine equations \cite[Proposition 3.3]{lubotzky:1988}. Alternately, the chain of isomorphisms \eqref{eq:LPSchain1} proves this fact by relating this Cayley graph to a quotient of a connected graph \cite[Theorem 7.4.3]{lubotzky:2010}: the infinite tree we shall describe in the next section.

The solutions $(x_0,x_1,x_2,x_3)$ and $(x_0,-x_1,-x_2,-x_3)$ correspond to elements of $S$ that are inverses in $\PSLt(\FF_p).$ Since $|S|=l+1$ this implies that the generators determine an undirected $(l+1)$-regular graph.

\subsection{Infinite tree of lattices}\label{subsubsect:infinitelatticetree}

Next we shall work over $\Ql.$ We give a description of the same graph in two ways: in terms of homothety classes of rank two lattices, and in terms of local double cosets of the multiplicative group of the Hamiltonian quaternion algebra. The description follows \cite[5.3, 7.4]{lubotzky:2010}. Let $B=B_{2,\infty}$ be the Hamiltonian quaternion algebra defined over $\QQ.$

First we review the construction of an $(l+1)$-regular infinite tree on homothety classes of rank two lattices in $\Ql\times \Ql $ following \cite[5.3]{lubotzky:2010}. The vertices of this infinite graph are in bijection with $\PGLt(\Ql)/\PGLt(\Zl).$ To talk about a finite graph, we shall then consider two subgroups $\Gamma(2)$ and $\Gamma(2p)$ in $B^{\times}/Z(B^{\times })$. It turns out that $\Gamma (2)$ acts simply transitively on the infinite tree, and orbits of $\Gamma (2p)$ on the tree are in bijection with the finite group $\Gamma (2)/\Gamma (2p).$ Under our assumptions the latter turns out to be in bijection with $\PSLt(\FF_p)$ above and the finite quotient of the tree is isomorphic to the Cayley graph above.

First we describe the infinite tree following \cite[5.3]{lubotzky:2010}. Consider the two dimensional vector space $\Ql\times \Ql$ with standard basis $\ve_1={}^t\langle 1,0 \rangle ,$ $\ve_2={}^t\langle 0,1 \rangle .$ A {\em{lattice}} is a rank two $\Zl$-submodule $L\subset \Ql\times \Ql.$ It is generated (as a $\Zl$-module) by two column vectors $\vu,\vv\in \Ql\times \Ql$ that are linearly independent over $\Ql.$ We shall consider homothety classes of lattices, i.e.\ we say lattices $L_1$ and $L_2$ are equivalent if there exists an $0\neq \alpha \in \Ql$ such that $\alpha L_1=L_2.$ Writing $\vu,\vv$ in the standard basis $\ve_1,\ve_2$ maps the lattice $L$ to an element $M_L\in\GLt(\Ql).$ Let $\vu_1,\vv_1,\vu_2,\vv_2\in \Ql\times \Ql$ and let $L_i=\Span_{\Zl}\{\vu_i,\vv_i\}$ ($i=1,2$) be the lattices generated by these respective pairs of vectors, with $M_{L_1}$ and $M_{L_2}$ the corresponding matrices. Let $M\in \GLt(\Ql)$ so that $M_{L_1}M=M_{L_2}.$ Then $L_1=L_2$ (as subsets of $\Ql\times \Ql$) if and only if $M\in \GLt(\Zl).$ It follows that the homothety classes of lattices are in bijection with $\PGLt(\Ql)/\PGLt(\Zl).$ Equivalently, we may say that $\PGLt(\Ql)/\PGLt(\Zl)$ acts simply transitively on homothety classes of lattices.

The vertices of the infinite graph $T$ are homothety classes of lattices. The classes $[L_1],[L_2]$ are adjacent in $T$ if and only if there are representatives $L_i'\in [L_i]$ ($i=1,2$) such that $L_2'\subset L_1'$ and $[L_1':L_2']=l.$ We show that this relation defines an undirected $(l+1)$-regular graph. By the transitive action of $\GLt(\Ql)$ on lattices we may assume that $L_1'=\Zl\times \Zl=\Span_{\Zl}\{\ve_1,\ve_2\},$ the {\em{standard lattice}} and $L_2'\subset \Zl\times \Zl.$ The map $\Zl\rightarrow\Zl/l\Zl\cong \FF_l$ induces a map from $\Zl\times \Zl$ to $\FF_l^2.$ Since the index of $L_2'$ in $\Zl\times \Zl$ is $l,$ the image of $L_2'$ is a one-dimensional vector subspace of $\FF_l^2.$ This implies that $L_2'\supset \{l\ve_1,l\ve_2\},$ i.e.\ $L_2'\supset lL_1'$ and the graph is undirected.\footnote{I.e.\ the adjacency relation defined above is symmetric.} Furthermore, since there are $l+1$ one-dimensional subspaces of $\FF_l^2,$ the graph is $(l+1)$-regular.

The $l+1$ neighbors of the standard lattice can be described explicitly by the following matrices:
\begin{equation}\label{eq:mtces_lattice_neighbours}
M_l=\tbtmtx{1}{0}{0}{l},\ M_h=\tbtmtx{l}{h}{0}{1}\text{ for }0\leq h\leq l-1
\end{equation}
For any of the matrices $M_t$ ($0\leq t\leq l$) the columns of $M_t$ span a different one-dimensional subspace of $\FF_l\times \FF_l.$ The matrices determine the neighbors of any other lattice by a  change of basis in $\Ql\times \Ql.$

By the above we can already see that $T$ is isomorphic to the graph on $\PGLt(\Ql)/\PGLt(\Zl)$ with edges corresponding to multiplication by generators \eqref{eq:mtces_lattice_neighbours} above. To show that $T$ is a tree it suffices to show that there is exactly one path from the standard lattice $\Zl\times \Zl $ to any other homothety class. This follows from the uniqueness of the Jordan--H\"{o}lder series in a finite cyclic $l$-group as in \cite[p.\ 69]{lubotzky:2010}.

In the next section, we show that the above infinite tree is isomorphic to a Cayley graph of a subgroup of $B^{\times }/Z(B^{\times }).$ In Section \ref{app:sect:explicitisom} we give an explicit bijection between the Cayley generators and the matrices given in \eqref{eq:mtces_lattice_neighbours} above.

\subsection{Hamiltonian quaternions over a local field}\label{subsubsect:localdoublecosets}

\newcommand{\vi}{{\mathbf{i}}}
\newcommand{\vj}{{\mathbf{j}}}
\newcommand{\vk}{{\mathbf{k}}}

To turn the above infinite tree into a finite, $(l+1)$-regular graph we shall define a group action on its vertices. Let $B$ be the algebra of Hamiltonian quaternions defined over $\QQ.$ Let $G'$ be the $\QQ$-algebraic group $B^{\times}/Z(B^{\times }).$ In this subsection we shall follow \cite[7.4]{lubotzky:2010} to define normal subgroups $\Gamma (2p)\subset \Gamma (2)$ of $\Gamma =G'(\ZZ[l^{-1}])$ such that $\Gamma (2)$ acts simply transitively on the graph $T.$ The quotient $\Gamma(2p)\backslash T $ will be isomorphic to the Cayley graph of the finite quotient group $\Gamma (2)/\Gamma(2p).$ This graph is isomorphic to the Cayley graph of $\PSLt(\FF_p)$ defined in Section \ref{subsubsect:CaylePres} above. Thus we have the following equation.
\begin{equation}\label{eq:chainofcongs_LPS_inHamiltsubsubsect}
\PSLt(\FF_p)\cong \Gamma(2p)\backslash \Gamma(2)\cong \Gamma(2p)\backslash T\cong \Gamma(2p)\backslash \PGLt(\Ql)/\PGLt(\Zl).
\end{equation}

We first define the groups $\Gamma ,$ $\Gamma (2),$ $\Gamma (2p)$ and then examine their relationship with $T.$ Recall that $B=B_{2,\infty }$, i.e.\ $B$ is ramified at $2$ and $\infty .$ For a commutative ring $R$ define $B(R)=\Span_{R}\{1,\vi, \vj,\vk\}$ where $\vi^2=\vj^2=-1$ and $\vi\vj=-\vj\vi=\vk.$ We introduce the notation $b_{x_0,x_1,x_2,x_3}:=x_0+x_1\vi+x_2\vj+x_3\vk.$ Recall that for $b=b_{x_0,x_1,x_2,x_3}$ we may define $\bar{b}=b_{x_0,-x_1,-x_2,-x_3}$ and the {\em{reduced norm}} of $b$ as $N(b)=b\bar{b}=x_0^2+x_1^2+x_2^2+x_3^2.$ For a (commutative, unital) ring $R$ an element $b\in B(R)$ is invertible in $B(R)$ if and only if $N(b)$ is invertible in $R.$ (Then $b^{-1}=(N(b))^{-1}\bar{b}.$) Furthermore
\begin{equation}\label{eq:commutatorb}
[b_{x_0,x_1,x_2,x_3},b_{y_0,y_1,y_2,y_3}]= 2(x_2y_3-x_3y_2)\vi + 2(x_3y_1-x_1y_3)\vj +2(x_1y_2-x_2y_1)\vk,
\end{equation}
and hence if $R$ has no zero divisors then $Z(B(R))=R.$  
In particular $Z(B^{\times }(\ZZ[l^{-1}]))=\{\pm l^k\mid k\in \ZZ\}.$

Recall that $S$ was the set of $l+1$ integer solutions of \eqref{eq:lsquaresumconds}. Any solution $x_0,x_1,x_2,x_3$ determines a $b=b_{x_0,x_1,x_2,x_3}\in B(\ZZ[l^{-1}])$ such that $N(b)=l .$ Since $l$ is invertible in $\ZZ[l^{-1}]$ we in fact have $b\in B^{\times }(\ZZ[l^{-1}]).$ Let $\Gamma =G'(\ZZ[l^{-1}])=B^{\times }(\ZZ[l^{-1}])/Z(B^{\times }(\ZZ[l^{-1}])) $ and let us denote the image of $S$ in $\Gamma $ by $S$ as well. Since $B^{\times }(\ZZ[l^{-1}])=\{b\in B(\ZZ[l^{-1}])\mid N(b)=l^k,\ k\in \ZZ \},$ if $[b]\in \Gamma $ for $b\in B^{\times }(\ZZ[l^{-1}])$ then it follows from \cite[Corollary 2.1.10]{lubotzky:2010} that $b$ is a unit multiple of an element of $\langle S \rangle .$ It follows that
$\Gamma =\langle S \rangle \{[1],[\vi],[\vj],[\vk] \}$ and the index of $\langle S \rangle$ in $\Gamma $ is $4.$ In fact observe that if $b\in S$ then $b^{-1}\in S$ and \cite[Corollary 2.1.11]{lubotzky:2010} states that $\langle S \rangle $ is a free group on $\frac{l+1}{2}$ generators. We shall see that $\langle S \rangle $ agrees with a congruence subgroup $\Gamma (2).$

Now let $N=2M$ be coprime to $l$ and let $R=\ZZ[l^{-1}]/N\ZZ[l^{-1}].$ The quotient map $\ZZ[l^{-1}]\rightarrow R$ determines a map $B(\ZZ[l^{-1}])\rightarrow B(R).$ This restricts to a map $B^{\times }(\ZZ[l^{-1}])\rightarrow B^{\times }(R).$ Observe that if $M=1$ then $B^{\times }(R)$ is commutative. If $M=p$ then the subgroup
$$Z:=\left\lbrace b_{x_0,0,0,0}\in B^{\times }(\ZZ[l^{-1}]/2p\ZZ[l^{-1}])\mid p\nmid x_0,2\nmid x_0 \right\rbrace$$
(cf.\ \cite[p. 266]{lubotzky:1988}) is central in $B^{\times }(R).$ Consider the commutative diagram:
\begin{equation}\label{eq:commdiag}
\begin{array}{ccccc}
B(\ZZ[l^{-1}])^{\times }&\longrightarrow & B^{\times }(\ZZ[l^{-1}]/2\ZZ[l^{-1}])& \longrightarrow & B^{\times }(\ZZ[l^{-1}]/2p\ZZ[l^{-1}])\\
\downarrow &&\downarrow&&\downarrow\\
\Gamma & \stackrel{\pi_{2} }{\longrightarrow} & B^{\times}(\ZZ[l^{-1}]/2\ZZ[l^{-1}]) &
\stackrel{\pi_{p} }{\longrightarrow} & B^{\times}(\ZZ[l^{-1}]/2p\ZZ[l^{-1}])/Z
\end{array}
\end{equation}
and define\footnote{The definition here agrees with the choices in \cite{lubotzky:1988} as well as $\Gamma (N)=\ker (G'(\ZZ[l^{-1}])\rightarrow G'(\ZZ[l^{-1}]/N\ZZ[l^{-1}]))$ in \cite{lubotzky:2010}. Here $G'=B^{\times }/Z(B^{\times })$ as a $\QQ$-algebraic group. Note however that by \eqref{eq:commutatorb} the center $Z(B^{\times }(R))$ for $R=\ZZ[l^{-1}]/N\ZZ[l^{-1}],$ $N=2M$ may not be spanned by $1+N\ZZ[l^{-1}].$ In fact from \eqref{eq:commutatorb} $B^{\times }(R)$ is commutative for $M=1$ and for $M=p$ we have $Z(B^{\times }(R))=Z\oplus [p]\vi +[p]\vj+[p]\vk .$ However the image of $\langle S\rangle $ in $ B^{\times }(R)$ is trivial if $M=1$ and intersects the center in $Z$ when $M=p.$} $\pi _{2p}:=\pi _{p}\circ \pi _{2}$ and $\Gamma (2):=\ker \pi _2$ and $\Gamma (2p)=\ker \pi _{2p}.$ Observe that by the congruence conditions (cf. \eqref{eq:lsquaresumconds}) $S\subseteq \Gamma $ is contained in $\Gamma (2)$ and in fact $\langle S\rangle =\Gamma (2)\supseteq\Gamma (2p).$ As mentioned above this implies that $\Gamma (2)$ is a free group with $\frac{l+1}{2}$ generators.

To see the action of $\Gamma (2)$ on $T$ note that $B$ splits over $\Ql$ and hence $B(\Ql)\cong M_2(\Ql).$ Since $-1\in (\FF_l^{\times})^2$ there exists an $\epsilon\in \Zl$ such that $\epsilon^2=-1.$ Then we have an isomorphism $\sigma : B(\Ql)\rightarrow M_2(\Ql)$ \cite[p. 95]{lubotzky:2010} given by
\begin{equation}\label{eq:unramlocalisomdef}
\sigma (x_0+x_1\vi+x_2\vj+x_3\vk) = \tbtmtx{x_0+x_1\epsilon}{x_2+x_3\epsilon}{-x_2+x_3\epsilon}{x_0-x_1\epsilon}.
\end{equation}
Observe that $\sigma (B^{\times }(\ZZ[l^{-1}]))\subseteq \GLt(\Ql)$ and $\sigma $ maps elements of the center into scalar matrices, and hence this defines an action of $\Gamma $ (and hence $\Gamma (2),\Gamma(2p)$) on $T.$ This action preserves the graph structure. Then we have the following.
Observe that $\sigma $ maps the elements of $\langle S \rangle\subseteq \Gamma $ into the congruence subgroup of $\PGLt(\Zl)$ modulo $2.$

\begin{proposition}\label{prop:simpletrans_on_tree}
\cite[Lemma 7.4.1]{lubotzky:2010} The action of $\Gamma (2)$ on the tree $T=\PGLt(\Ql)/\PGLt(\Zl)$ is simply transitive (and respects the graph structure).
\end{proposition}
\begin{proof}
See {\em{loc.cit.}}\ for details of the proof. Transitivity follows from the fact that $T$ is connected and elements of $S$ map a vertex of $T$ to its distinct neighbors. The group $\Gamma (2)=\langle S\rangle$ is a discrete free group, hence its intersection with a compact stabilizer $\PGLt(\Zl)$ is trivial. This implies that the neighbors are distinct and the stabilizer of any vertex is trivial.
\end{proof}

The above implies that the orbits of $\Gamma (2p)$ on $T$ have the structure of the Cayley graph $\Gamma (2)/\Gamma (2p)$ with respect to the generators $S.$ We can see from the maps in \eqref{eq:commdiag} that $\Gamma (2)/\Gamma (2p)$ is isomorphic to a subgroup of $G'(\ZZ/2p\ZZ)\cong G'(\ZZ/2\ZZ)\times G'(\ZZ/p\ZZ).$ (This last isomorphism follows from the Chinese Remainder Theorem.) Since the image of $\Gamma (2)$ in $G'(\ZZ/2\ZZ)$ is trivial, we may identify $\Gamma (2)/\Gamma (2p)$ with a subgroup of $G'(\ZZ/p\ZZ).$ Here $G'(\ZZ/p\ZZ)\cong \PGLt(\FF_p).$ (For an explicit isomorphism take an analogue of $\sigma $ in \eqref{eq:unramlocalisomdef} with $\epsilon \in \ZZ/p\ZZ$ such that $\epsilon ^2=-1.$) The image of $\Gamma (2)$ agrees with $\PSLt(\FF_p)$ as a consequence of the Strong Approximation Theorem \cite[Lemma 7.4.2]{lubotzky:2010}. We shall discuss this in the next section.

We summarize the contents of this section.

\begin{theorem}\label{thm:LPS-summary}
\cite[Theorem 7.4.3]{lubotzky:2010}
Let $l$ and $p$ be primes so that $l\equiv p\equiv 1\mod{4}$ and $l$ is a quadratic residue modulo $2p.$ Let $S\subset \PSLt(\FF_p)$ be the $(l+1)$-element set corresponding to the solutions of \eqref{eq:lsquaresumconds} via the map \eqref{eq:elementsofSmappedtomtcesmodp} and $Cay(\PSLt(\FF_p),S)$ the Cayley graph determined by the set of generators $S$ on the group $\PSLt(\FF_p).$ Let $T$ be the graph on $\PGLt(\Ql)/\PGLt(\Zl)$ with edges corresponding to multiplication by elements listed in \eqref{eq:mtces_lattice_neighbours}. Let $B$ be the Hamiltonian quaternion algebra over $\QQ$ and $\Gamma (2p)$ the kernel of the map $\pi_{2p}$ in \eqref{eq:commdiag} (a cocompact congruence subgroup). Then $\Gamma (2p)$ acts on the infinite tree $T$ and we have the following isomorphism of graphs:
\begin{equation}\label{eq:thm:LPS}
Cay(\PSLt(\FF_p),S)\cong \Gamma (2p) \backslash \PGLt(\Ql)/\PGLt(\Zl).
\end{equation}
These are connected, $(l+1)$ regular, non-bipartite, simple, graphs on $\frac{p^3-p}{2}$ vertices.
\end{theorem}

\subsection{Explicit isomorphism between generating sets}\label{app:sect:explicitisom}

We have seen above that the LPS graph can be interpreted as a finite quotient of the infinite tree of homothety classes of lattices. In this case, the edges are given by matrices that take a $\Zl$-basis of one lattice to a $\Zl$-basis of one of its neighbors. On the other hand, the edges can be given in terms of the set of generators $S.$
Proposition \ref{prop:simpletrans_on_tree} states that $\langle \sigma(S)\rangle =\Gamma (2)\subset G'(\ZZ[l^{-1}])$ acts simply transitively on the tree $T.$ The proof of the proposition (cf. \cite[Lemma 7.4.1]{lubotzky:2010}) implicitly shows that there exists a bijection between elements of $\sigma(S)\subset \PGLt(\Zl)$ and the matrices given in \eqref{eq:mtces_lattice_neighbours}.

In this section we wish to make this bijection more explicit. For a fixed $\alpha \in S$ we find the matrix from the list \eqref{eq:mtces_lattice_neighbours} determining the same edge of $T$. As in Section \ref{subsubsect:localdoublecosets} we write $\sigma (\alpha )\in \PGLt(\Zl)$ for the elements of $\sigma(S).$
This amounts to finding the matrix $M$ from the list in \eqref{eq:mtces_lattice_neighbours} such that $\sigma(\alpha) ^{-1}M\in \PGLt(\Zl).$

To pair up matrices from \eqref{eq:mtces_lattice_neighbours} with the corresponding elements of $S,$ we introduce the following notation. Let us number the solutions to $\alpha \overline{\alpha }=l$ as $\alpha _0,\ldots ,\alpha_{l-1},\alpha_l$ so that we have the correspondence $\sigma(\alpha _h) ^{-1}M_h\in \PGLt(\Zl)$ for $0\leq h\leq l.$
By giving an explicit correspondence, we mean that given an $\alpha \in \sigma ^{-1}(S),$ we determine $0\leq h\leq l$ such that $\alpha =\alpha _h.$

Elements of $\sigma(S)\subset \PGLt(\Zl)$ are given in terms of an $\epsilon \in \Zl$ such that $\epsilon^2=-1.$ Let $a,$ $b$ be the positive integers such that $a^2+b^2=l$ and $a$ is odd. Let $0\leq e\leq l-1$ so that $eb=a.$ Then in $\Zl$ we have either $\epsilon\in e+l\Zl$ and $\epsilon^{-1}=-\epsilon\in -e+l\Zl$ or $\epsilon\in -e+l\Zl$ and $\epsilon^{-1}=-\epsilon\in e+l\Zl.$

Let $\alpha =x_0+x_1\vi+x_2\vj+x_3\vk$ so that $\sigma (\alpha )\in S,$ and $a,b,e, \epsilon$ are as above. Let
$$\alpha_h =x_0^{(h)}+x_1^{(h)}\vi+x_2^{(h)}\vj+x_3^{(h)}\vk$$ for $0\leq h\leq l.$  Here $x_0,x_1,x_2,x_3$ are integers; it is convenient to think about them (as well as $x_0^{(h)},x_1^{(h)},x_2^{(h)},x_3^{(h)}$ for $0\leq h\leq l$) as being in $\ZZ \subset \Zl .$
Then
\begin{equation}\label{eq:sigmainv_alpha}
\sigma(\alpha) ^{-1}=\frac{1}{l}\tbtmtx{x_0-x_1\epsilon}{-x_2-x_3\epsilon}{x_2-x_3\epsilon}{x_0+x_1\epsilon}
\end{equation}
and
\begin{equation}\label{eq:timesmtx2} 
\begin{split}
\sigma(\alpha) ^{-1}\cdot \tbtmtx{l}{h}{0}{1}= & \tbtmtx{x_0-x_1\epsilon}{l^{-1}\left(h(x_0-x_1\epsilon)+(-x_2-x_3\epsilon)\right)}{x_2-x_3\epsilon}{l^{-1}\left(h(x_2-x_3\epsilon)+(x_0+x_1\epsilon)\right)}\\
\sigma(\alpha) ^{-1}\cdot \tbtmtx{1}{0}{0}{l}=& \tbtmtx{l^{-1}(x_0-x_1\epsilon)}{-x_2-x_3\epsilon}{l^{-1}(x_2-x_3\epsilon)}{x_0+x_1\epsilon}
\end{split}
\end{equation}

Then by \eqref{eq:timesmtx2} we have that $x_0^{(l)}-x_1^{(l)}\epsilon$ and $x_2^{(l)}-x_3^{(l)}\epsilon$ are in $l\Zl.$ Hence $x_0^{(l)}\in x_1^{(l)}\epsilon+l\Zl,$ and thus $(x_0^{(l)})^2\in (x_1^{(l)}\epsilon)^2+l\Zl=-x_1^2+l\Zl,$ whence $(x_0^{(l)})^2+(x_1^{(l)})^2\in l\Zl .$ Note that since $(x_0^{(l)})^2+(x_1^{(l)})^2+(x_2^{(l)})^2+(x_3^{(l)})^2=l$ and $x_0$ is positive, this implies that $(x_0^{(l)})^2+(x_1^{(l)})^2=l$ and $(x_2^{(l)})^2+(x_3^{(l)})^2=0,$ i.e. $x_2^{(l)}=x_3^{(l)}=0$ and $x_0^{(l)}=a,$ $|x_1^{(l)}|=b.$ Note that by the assumptions in Section \ref{subsubsect:CaylePres}, $a\pm bi,a\pm bj,a\pm bk\in S.$ A straightforward computation now shows the following.
\begin{equation}\label{eq:corresp:explicitpart1}
\begin{split}
\epsilon\in e+l\Zl \Rightarrow \alpha_l=a+b\vi,\ \alpha _0=a-b\vi,\ \alpha _e=a-b\vj,\ \alpha _{l-e}=a+b\vj,\ \alpha_1=a-b\vk,\ \alpha_{l-1}=a+b\vk\\
\epsilon\in -e+l\Zl\Rightarrow \alpha_l=a-b\vi,\ \alpha _0=a+b\vi,\ \alpha _e=a-b\vj,\ \alpha _{l-e}=a+b\vj,\ \alpha_1=a+b\vk,\ \alpha_{l-1}=a-b\vk
\end{split}
\end{equation}

Now let us assume that for $\alpha = x_0+x_1\vi+x_2\vj+x_3\vk$ we have that $x_0-x_1\epsilon\notin l\Zl.$ This implies that
It remains to determine the $h$ such that $\alpha =\alpha _h$ when $\alpha $ is not one of the solutions covered by \eqref{eq:corresp:explicitpart1}. In that case, we may assume $h\notin\{0,1,e,l-e,l-1,l\}$ and we have
\begin{equation}\label{eq:genh_eq1}
h(x_0-x_1\epsilon)+(-x_2-x_3\epsilon)\in l\Zl;
\end{equation}
\begin{equation}\label{eq:genh_eq2}
h(x_2-x_3\epsilon)+(x_0+x_1\epsilon)\in l\Zl.
\end{equation}
A straightforward computation based on $\alpha\overline{\alpha}=l$ shows that \eqref{eq:genh_eq1} and \eqref{eq:genh_eq2} are satisfied by the same element in $\FF_l=\ZZ/l\ZZ.$ The element
\begin{equation}\label{eq:solutiontoh}
\overline{h}=\frac{x_2+x_3\epsilon}{x_0-x_1\epsilon} \in \FF_l
\end{equation}
is well defined, since $x_0-x_1\epsilon\notin l\Zl ,$ furthermore, it uniquely determines an $0\leq h\leq l.$ For a fixed $\alpha $ not covered by \eqref{eq:corresp:explicitpart1}, one may thus find $h$ such that $\alpha =\alpha _h.$

We give two explicit examples.

\begin{example}\label{eg:corresp_l5}
When $l=5,$ then $a=1,$ $b=2$ and $e=3.$ Then \eqref{eq:example_corresp_l5} gives the bijection between the list in \eqref{eq:mtces_lattice_neighbours} and solutions of $\alpha\overline{\alpha}=5$ in $B(\QQ_5).$ In this case the list in \eqref{eq:corresp:explicitpart1} is exhaustive.
\begin{equation}\label{eq:example_corresp_l5}
\begin{array}{l|c||c|c|c|c|c|c}
\multicolumn{2}{c||}{h}  & 0 & 1& 2& 3& 4& 5 \\ \hline
 \epsilon \in 3+5\ZZ_5 & \multirow{2}{*}{$\alpha_h$} & 1-2\vi & 1-2\vk & 1+2\vj & 1-2\vj & 1+2\vk & 1+2\vi\\ \hhline{-~------}
\epsilon \in 2+5\ZZ_5 & & 1+2\vi & 1+2\vk & 1+2\vj & 1-2\vj & 1-2\vk & 1-2\vi
\end{array}
\end{equation}
\end{example}

\begin{example}\label{eg:corresp_l13}
When $l=13,$ we have $a=3,$ $b=2$ and $e=8.$ The cases listed in \eqref{eq:corresp:explicitpart1} are no longer exhaustive.
The correspondence is given in Table \ref{table:example_corresp_l13_8_and_5}. 
\begin{table}[h!]
\begin{center}
$\displaystyle \begin{array}{c|c}
h & \alpha _h \\ \hline \hline
0& 3-2\vi\\ \hline
1& 3-2\vk\\ \hline
2& 1 -2\vi -2\vj -2\vk\\ \hline
3& 1 -2\vi +2\vj -2\vk\\ \hline
4& 1 +2\vi +2\vj +2\vk\\ \hline
5& 3+2\vj\\ \hline
6& 1 +2\vi -2\vj +2\vk\\ \hline
7& 1 +2\vi +2\vj -2\vk\\ \hline
8& 3-2\vj\\ \hline
9& 1 +2\vi -2\vj -2\vk\\ \hline
10& 1 -2\vi -2\vj +2\vk\\ \hline
11& 1 -2\vi +2\vj +2\vk\\ \hline
12& 3+2\vk\\ \hline
13 & 3+2\vi
\end{array}$ \hfil $\displaystyle \begin{array}{c|c}
h & \alpha _h \\ \hline \hline
0& 3+2\vi\\ \hline
1& 3+2\vk\\ \hline
2& 1 +2\vi -2\vj +2\vk\\ \hline
3& 1 +2\vi +2\vj +2\vk\\ \hline
4& 1 -2\vi +2\vj -2\vk\\ \hline
5& 3+2\vj\\ \hline
6& 1 -2\vi -2\vj -2\vk\\ \hline
7& 1 -2\vi +2\vj +2\vk\\ \hline
8& 3-2\vj\\ \hline
9& 1 -2\vi -2\vj +2\vk\\ \hline
10& 1 +2\vi -2\vj -2\vk\\ \hline
11& 1 +2\vi +2\vj -2\vk\\ \hline
12& 3-2\vk\\ \hline
13 & 3-2\vi
\end{array}$
\end{center}
\caption{The correspondence when $\epsilon \in 8+13\ZZ_{13}$ (left) and when $\epsilon \in 5+13\ZZ_{13}$ (right).}
\label{table:example_corresp_l13_8_and_5}
\end{table}

\end{example}

\section{Strong Approximation}\label{subsect:StrongApprox}
In this section we briefly explain the significance of Strong Approximation to Ramanujan graphs and particularly the LPS graphs above. As discussed in Section \ref{subset:autoprelim} we may consider $G(\A),$ the adelic points of a linear algebraic group $G$ defined over $\QQ.$
The group $G(\QQ)$ embeds diagonally into $G(\A),$ and it is a discrete subgroup. The groups $G(\QQ_v)$ are also subgroups of $G(\A),$ and $G(\A)$ has a well-defined projection onto $G(\QQ_v).$ Similarly, for a finite set of places $S$ we may take $G_S,$ the direct product of $G(\QQ_v)$ for $v\in S.$

Strong Approximation (when it holds) is the statement that for a group $G$ and a finite set of places $S$ the subgroup $G(\QQ)G_S$ is {\em{dense}} in $G(\A).$ This implies that
\begin{equation}\label{eq:SA_anyK}
G(\A)=G(\QQ)G_SK\text{ for any open subgroup }K\leq G(\A).
\end{equation}

For example, Strong Approximation holds for $G=\SLt$ and any set of places $S=\{v\}.$ However, in the form written above it {\em{does not hold}} for $\GLt$ or $\PGLt .$ However one can prove results similar to \eqref{eq:SA_anyK} for $\GLt$ adding restrictions on the subgroup $K$:
\begin{equation}\label{eq:SA_notanyK}
G(\A)=G(\QQ)G_SK\text{ for an open subgroup }K\leq G(\A)\text{ if $K$ is ``sufficiently large.''}
\end{equation}
Here we shall have
\begin{equation}
K=\prod_{v\notin S}K_v;\ K_v\leq G(\ZZ_v)
\end{equation}
and the condition of being ``sufficiently large'' can be made precise by requiring that the determinant map $\det: K_v\rightarrow \ZZ_v^{\times }$ be surjective for all $v \not \in S$.

Strong Approximation holds for the algebraic group of elements of a quaternion algebra of unit norm \cite[Th\'{e}or\`{e}me 4.3]{vigneras2006arithmetique}. We shall use this statement to prove a statement like \eqref{eq:SA_notanyK} for the algebraic group of invertible quaternions. A similar statement then holds for $G'=B^{\times}/Z(B^{\times })$ and a subgroup $K'$ that is not quite ``large enough.''
 The implications for Pizer graphs and LPS graphs will be discussed in Sections \ref{subsubsect:SA_LPS} and \ref{subsubsect:SA_Pizer} below.

These statements coming from Strong Approximation are crucial for proving that the various constructions produce Ramanujan graphs. 
As seen in Section \ref{subset:autoprelim} the Ramanujan property of a graph can be expressed in terms of its eigenvalues. Given a graph (constructed e.g.\ via local double cosets as seen above) the Strong Approximation theorem can be used to relate its spectrum to the representation theory of $G(\A).$ In that context a theorem of Deligne resolves the issue by proving a special case of the Ramanujan conjecture (see \cite[Theorem 6.1.2, Theorem A.1.2, Theorem A.2.14]{lubotzky:2010} and \cite{deligne11formes}).

\subsection{Approximation for invertible quaternions}\label{subsubsect:SA_Bcross}

The argument below is adapted from \cite[Section 3]{gelbart1975automorphic} and \cite[6.3]{lubotzky:2010}.\footnote{In fact, since at every split place $v$ we have $B^{\times }(\Qv) \cong \GLt(\Qv)$ with the reduced norm on $B^{\times }$ corresponding to the determinant on $\GLt$ \cite[p.\ 3]{vigneras2006arithmetique} this is the ``same argument at all but finitely many places.''}

Let $B$ be a (definite) quaternion algebra over $\QQ$, $B^{\times }$ its invertible elements and $B^1=\{b\in B\mid N(b)=1\}$ its elements of reduced norm $1,$ recall $N(b)=b\bar{b}.$ Let $l$ be a prime where $B$ is split. Then by \cite[Th\'eor\`{e}me 4.3]{vigneras2006arithmetique} we have that $B^1(\QQ)B^1(\Ql)$ is dense in $B^1(\A)$ thus $B^1(\A)=B^1(\QQ)B^1(\Ql)K$ for any open subgroup $K\leq B^1(\A ).$ An open subgroup $K\leq B^1(\A)$ is of the form $K=\prod_{v} K_v$ where $K_v\leq B^1_v$ is open and $K_v=B^1(\Zv)$ for all but finitely many places $v.$ It follows that given {\em{any}} open subgroups $K_v^{(B^1)}\leq B^1(\Zv)$ ($v\neq l$) such that $K_v^{(B^1)}=B^1(\Zv)$ for all but finitely many places $v$ we have that
\begin{equation}\label{eq:SA_H_genK}
B^1(\A )=B^1(\QQ)B^1(\Ql)\prod_{v\neq l} K_v^{(B^1)}.
\end{equation}

To make a similar statement for $B^{\times }$ it will be necessary to impose a restriction on the open subgroups $K_v.$ 
\begin{theorem}\label{thm:SA_Bcross}
Let $K_v\leq B^{\times }(\Zv)$ for every place $l\neq v<\infty$ so that $K_v=B^{\times }(\Zv)$ for all but finitely many $v,$ and the norm map $N:K_v\rightarrow \Zv^{\times }$ is surjective for every place $v$. Then
\begin{equation}\label{eq:weakerthan_SA_forBcross}
B^{\times }(\A)=B^{\times }(\QQ)B^{\times }(\RR) B^{\times }(\Ql)\prod_{l\neq v<\infty}K_v.
\end{equation}
\end{theorem}
 Note that by \cite[Lemma 13.4.6]{voight:quatbook} the norm map $N: B^{\times }(\Zv)\rightarrow \Zv^{\times }$ is surjective for every nonarchimedean $v.$
\begin{proof}
Let $b\in B^{\times }(\A),$ we need to show $b$ is contained on the right-hand side.
To write $b$ as a product according to the right-hand side of \eqref{eq:weakerthan_SA_forBcross} we shall use \eqref{eq:SA_H_genK}, strong approximation for $B^1.$ Observe first that it suffices to show that any $b\in B^{\times}(\A)$ can be written as
\begin{equation}\label{eq:suffices_weakappr_forb}
b=rhk,\text{ where }r\in B^{\times}(\QQ),\ h\in B^1(\A),\text{ and }k\in B^{\times }(\RR) B^{\times }(\Ql)\prod_{l\neq v<\infty}K_v.
\end{equation}
This is because the intersections $K_v\cap B^1(\Qv)$ are open subgroups of $B^1(\Zv)$ (and $B^{\times }(\Zv)\cap B^1(\Zv)=B^1(\Zv)$ at all but finitely many places). It thus follows from \eqref{eq:SA_H_genK} (choosing $K_v^{(B^1)}:=K_v\cap B^1(\Qv)$) that the factor $h\in B^1(\A)\subseteq B^{\times }(\A)$ from \eqref{eq:suffices_weakappr_forb} is contained on the right-hand side of \eqref{eq:weakerthan_SA_forBcross}. It follows that then $b=rhk$ is contained on the right-hand side of \eqref{eq:weakerthan_SA_forBcross} as well. (Note that here the factors of $h$ and $k$ belonging to different components $B^{\times }(\Qv)$ commute.)

So we must show that any $b\in B^{\times }(\A)$ decomposes as in \eqref{eq:suffices_weakappr_forb}. Let $b=(b_v)_v$ for $b_v\in B^{\times }(\Qv )$ and set $n_v:=N(b_v).$ For all but finitely many places $v$ we have $b_v\in  B^{\times }(\Zv )$ and hence $n_v\in \Zv^{\times }.$ At a finite set $T$ of finite places we may write $n_v\in v^{m_v}\Zv^{\times }.$ Let us take
\begin{equation}\label{eq:choicenQ}
n_{\QQ}=\prod _{v\in T} v^{m_v}.
\end{equation}
Then $n_{\QQ}\in \QQ_{>0},$ $n_{\QQ}\in \Zv^{\times }$ for every $v\notin T,$ $v<\infty$ and hence $n_{\QQ}^{-1}n_v\in \Zv^{\times }$ for {\em{every}} finite place $v.$

It is a fact that there is an $r\in B^{\times }(\QQ)$ such that $N(r)=n_{\QQ}.$ Then for this $r$ we have that the norm of $r^{-1}b\in B^{\times }(\A)$ is in $\Zv^{\times }$ for every finite place $v.$

Let us write $(r^{-1}b)_v$ for the component of $r^{-1}b\in B^{\times }(\A)$ at a place $v.$ There exists a $k\in B^{\times }(\RR) B^{\times }(\Ql)\prod_{l\neq v<\infty}K_v,$ $k=(k_v)_v$ such that $k_l=(r^{-1}b)_l$ and $k_{\infty }=(r^{-1}b)_{\infty }$ and $N(k_v)=N((r^{-1}b)_v)$ every other place. This follows from the fact that the norm map $N:K_v\rightarrow \Zv^{\times }$ is surjective.

Now let $h=r^{-1}bk^{-1}.$ We show $h\in B^1(\A ).$ Write $h=(h_v)_v$ for $h_v\in B^{\times }(\Qv).$ It follows from the choice of $k$ that $h_l$ and $h_{\infty }$ are the identity element of $B^{\times}(\Ql)$ and $B^{\times}(\RR)$ respectively, and $N(h_v)=1$ at every other place $v.$ This implies that indeed $h\in B^1(\A ).$ This completes the proof that a decomposition as in \eqref{eq:suffices_weakappr_forb} exists, and in turn the proof of \eqref{eq:weakerthan_SA_forBcross}.
\end{proof}

\subsection{Strong Approximation for LPS graphs}\label{subsubsect:SA_LPS}

This section is based on \cite[6.3]{lubotzky:2010}. (In particular, we recall and elaborate on the proof of the first statements in \cite[Proposition 6.3.3]{lubotzky:2010} in the special case when $N=2p.$ This is relevant to understanding the last step in \eqref{eq:LPSchain1}.) We apply a similar formula to \eqref{eq:weakerthan_SA_forBcross} with a particular choice of open subgroups $K_v'$ to prove a statement that relates double cosets such as in \eqref{eq:chainofcongs_LPS_inHamiltsubsubsect} to adelic double cosets. Let $B=B_{2,\infty}$ be the algebra of Hamiltonian quaternions, ramified at $2$ and $\infty .$ Recall from Section \ref{subsubsect:localdoublecosets} that $G'$ is the $\QQ$-algebraic group $B^{\times}/Z(B^{\times }).$
Let us fix the prime $l\equiv 1\mod{4}$ as in Section \ref{subsect:LPS}.
In a similar manner to the proof of \eqref{eq:weakerthan_SA_forBcross} is follows that
\begin{equation}\label{eq:G'A-Approx}
G'(\A)=G'(\QQ)G'(\RR) G'(\Ql)\prod_{l\neq v<\infty} G'(\Zv).
\end{equation}

Recall that since $B$ splits at $l$ we have $G'(\Ql)\cong \PGLt(\Ql).$ We wish to have a statement similar to \eqref{eq:G'A-Approx} above, replacing $G'(\Zv)$ at $v=2$ and $v=p$ by congruence subgroups $K_2'$ and $K_p'.$ (This $p$ is the one fixed above in Section \ref{subsect:LPS}.) Then isomorphism will no longer hold, but the right-hand side will be a finite index normal subgroup of $G'(\A).$

The choice of the smaller subgroups $K_2'$ and $K_p'$ is as follows. For $v\in \{2,p\}$ let
\begin{equation}\label{eq:compactsgdef_LPS_2p}
K_v'=\ker \left( G'(\Zv)\rightarrow G'(\Zv/v\Zv) \right).
\end{equation}
Here $\Zv/v\Zv=\FF_v$ is a finite field, hence $G'(\Zv/v\Zv)$ is finite. It follows that the index $[K_v:K_v']$ is finite. In fact since $B_{2,\infty}$ splits over $p$ we have that $G'(\Zp/v\Zp)\cong \PGLt(\FF_p),$ hence $[K_p:K_p']=p(p^2-1).$ At $v=2$ we have $G'(\FF_2)=B^{\times }(\FF_2)$ hence $[K_2:K_2']=8.$

Let us set $K_v'$ as above if $v\in \{2,p\}$ and $K_v'=K_v=G'(\Zv)$ otherwise, and let us define
\begin{equation}\label{eq:2padelicgroupdef}
H_{2p}:=\left(G'(\QQ)G'(\RR) G'(\Ql)\prod_{l\neq v<\infty} K_v' \right).
\end{equation}
By \cite[Proposition 6.3.3]{lubotzky:2010}  Strong Approximation proves that $H_{2p}$ is a finite index normal subgroup of $G'(\A)$.

From the definition of $H_{2p}$ in equation \eqref{eq:2padelicgroupdef} we have a surjection from
\[
G'(\QQ_l)\rightarrow G'(\QQ)\backslash H_{2p}/G'(\RR)\prod_{l\neq v<\infty} K_v'.
\]
If $g_l$ and $g_l' \in G'(\QQ_l)$ are mapped to the same coset on the right hand side then there exists $g_q\in G'(\QQ), g_r \in G'(\RR)$ and $k=\prod_{l\neq v < \infty} k_v\in \prod_{l\neq v<\infty} K_v'$ such that $g_l=g_q g_l' g_r k$. This is equivalent to saying $g_l=g_q g_l'$ and $g_q \in K_v'$ for all $l\neq v<\infty$. By the definitions of the $K_v'$s this last condition implies $g_q \in \Gamma(2p)$. Thus we see that
\begin{equation}\label{eq:localglobalLPS}
\Gamma(2p)\backslash G'(\QQ_l)/G'(\ZZ_l) \cong G'(\QQ)\backslash H_{2p}/G'(\RR)\prod_{v<\infty} K_v'.
\end{equation}

Strong approximation in the manner discussed above is used to prove that  LPS graphs are Ramanujan. First one shows that the finite $(l+1)$-regular graph $\Gamma(2p)\backslash T$ is Ramanujan if and only if all  irreducible infinite-dimensional  unramified unitary representations of $\PGL_2(\QQ_l)$ that appear in $L^2(PGL_2(\QQ_l)/\Gamma(2p))$
are tempered
\cite[Corollary 5.5.3]{lubotzky:2010}. Then by the isomorphism above which follows from Strong Approximation, one can extend a representation $\rho'_l$ of $\PGL_2(\QQ_l)$ to an automorphic representation $\rho'$ of $G'(\A)$ in $L^2(G'(\QQ)\backslash G'(\A))$. By the Jacquet--Langlands correspondence, $\rho'$ corresponds to a cuspidal representation $\rho$ of $\PGL_2(\A)$ in $L^2(\PGL_2(\QQ)\backslash \PGL_2(\A))$ such that $\rho_v$ is discrete series for all $v$ where $B$ ramifies (so in our case, $2$ and $\infty$) \cite[Theorem 6.2.1]{lubotzky:2010}. Finally, Deligne has proved the Ramanujan--Peterson conjecture in this case of holomorphic modular forms \cite[Theorem 6.1.2]{lubotzky:2010}, \cite{deligne11formes}, \cite{deligne1974conjecture}
which says that for $\rho$ a cuspidal representation of $\PGL_2(\A)$ in $L^2(\PGL_2(\QQ)\backslash \PGL_2(\A))$ with $\rho_\infty$ discrete series, $\rho_l$ is tempered \cite[Theorems 7.1.1 and 7.3.1]{lubotzky:2010}.
Under the Jacquet--Langlands correspondence, the adjacency matrix of  our graph $X$ corresponds to the Hecke operator $T_l$ \cite[5.3]{lubotzky:2010} and the Ramanujan conjecture is equivalent to saying that $|\lambda|\leq 2\sqrt{l}$ for all of its eigenvalues $\lambda\neq \pm(l+1) $.

\subsection{Strong Approximation for Pizer graphs}\label{subsubsect:SA_Pizer}

Now we turn to discussing how strong approximation is useful in establishing the bijections in \eqref{eq:whatEyalwrote}. In Section \ref{sect:Pizer} we will discuss Pizer's construction of Ramanujan graphs. These graphs are isomorphic to supersingular isogeny graphs. Their vertex set is the class group of a maximal order $\O$ in the quaternion algebra $B_{p,\infty}.$ This set is in bijection with an adelic double coset space, which in turn is in bijection with a set of local double cosets.

Let $B=B_{p,\infty }$ be a quaternion algebra (over $\QQ$) ramified exactly at $\infty $ and at a finite prime $p.$ At every finite prime $v$, $B(\Qv)$ has a unique maximal order up to conjugation \cite[Lemme 1.4]{vigneras2006arithmetique}. Given a maximal order $\O$ of $B,$ one may define the adelic group $B^{\times }(\A _f)$ as a restricted direct product of the groups $B^{\times }(\Qv)$ over the finite places, with respect to $\O_v^\times.$ (Recall that this means that any element of $B^{\times }(\A _f)$ is a vector indexed by the finite places $v;$ the component at $v$ is in $B^{\times }(\Qv)$ and in fact in $\O_v^\times$ at all but finitely many places.) This adelic object does not in fact depend on the choice of the maximal ideal $\O.$ In particular, at any prime $l\neq p$ where $B$ splits we have $B^{\times}(\Ql)\cong \GLt(\Ql)$ and $\O_l^\times\cong\GLt(\Zl).$

Let us now fix a prime $l$ where $B$ splits. The same argument as in Section \ref{subsubsect:SA_Bcross} works restricted to $B^{\times }(\A_f)$ (the finite ad\`{e}les). It follows that we have
\begin{equation}\label{eq:SA_finadeles_B*}
B^{\times }(\A_f)=B^{\times }(\QQ)B^{\times }(\Ql)\prod_{l\neq v<\infty}B^{\times }(\Zv).
\end{equation}

\begin{proposition}\label{prop:SA_Pizer}
We have the bijections (cf.\ \cite[(1)]{lauter:2009r})
\begin{equation}\label{eq:Pizer_doublecoset-bijection}
\begin{split}
B^{\times }(\QQ)\backslash B^{\times }(\A_f)\slash \prod_{v<\infty}B^{\times }(\Zv)\cong & (\O(\ZZ[l^{-1}]))^{\times }\backslash B^{\times }(\Ql)\slash B^{\times }(\Zl)\\
\cong &  (\O(\ZZ[l^{-1}]))^{\times }\backslash \GLt(\Ql)\slash \GLt(\Zl).
\end{split}
\end{equation}
\end{proposition}
\begin{proof}
The first bijection follows from \eqref{eq:SA_finadeles_B*} and an argument similar to the proof of \eqref{eq:localglobalLPS}. Indeed, \eqref{eq:SA_finadeles_B*} implies that there is a surjection
\begin{equation}\label{eq:Pizer_SA_surjfirst}
B^{\times }(\Ql)\rightarrow B^{\times }(\QQ)\backslash B^{\times }(\A_f)\slash \prod_{l\neq v<\infty}B^{\times }(\Zv).
\end{equation}
Now two elements $g_l,g_l'\in B^{\times }(\Ql)$ land in the same double coset via this bijection if and only if $g_l=g_qg_l'k$ in $B^{\times }(\A_f).$ Then $g_l=g_qg_l'$ (from equality at the place $l$) and $g_q\in B^{\times }(\Zv)$ (from equality at the places $l\neq v<\infty$). Consider the element $g_q\in B(\QQ),$ for example in terms of its coordinates in the standard basis $\{1,\vi, \vj,\vk\}$ of $B.$ Since $g_q\in B^{\times }(\Zv)$ we have that $g_q\in \O(\ZZ[l^{-1}]),$ and $g_q\in B^{\times }(\Ql)$ implies that in fact $g_q\in (\O(\ZZ[l^{-1}]))^{\times }.$ This completes the proof of the first bijection in \eqref{eq:Pizer_doublecoset-bijection}.

Now the second bijection follows from the fact that $B$ splits at the prime $l$ and hence $B^{\times }(\Ql)\cong \GLt(\Ql)$ with the unique maximal order $\GLt(\Zl).$
\end{proof}

Finally, we wish to also address the bijection between the adelic double coset object and the class group of the maximal order $\O .$ This fact follows from the fact that ideals of $\O$ are locally principal. We omit defining ideals of an order $\O$ or defining the class group here and instead refer the reader to \cite[\S 4]{vigneras2006arithmetique}, \cite[\S 2.3]{chenevier} or \cite{voight:quatbook}. For the statement about the bijection between the class group $\Cl(\O)$ and the adelic double cosets in \eqref{eq:Pizer_doublecoset-bijection} above, see for example \cite[Theorem 2.6]{chenevier}.

\section{Pizer Graphs}\label{sect:Pizer}

In this section we give an overview of Pizer's \cite{pizer:1998} construction of a Ramanujan graph. The graphs constructed by Pizer are isomorphic to the graphs of supersingular elliptic curves over $\FF_{p^2}$ \cite[Section 2]{lauter:2009r}. These graphs were considered by Mestre \cite{mestre:1986} and Ihara \cite{ihara:1966} before (cf.\ \cite{jao2005all}), but Pizer's construction reveals their connection to quaternion algebras, proving their Ramanujan property. In Section \ref{subsect:LPS_Pizer_compare} we shall compare the resulting graphs to the LPS construction described above.

Pizer's description is in terms of a quaternion algebra and a pair of prime parameters $p,l.$ We shall aim to keep technical details to a minimum, and focus on the choice of quaternion algebra and parameters. This elucidates the connection with the LPS construction. Recall that the meaning of the parameters is similar in both cases: the resulting graphs are $(l+1)$-regular and their size depends on the value of $p.$ Varying $p$ (subject to some constraints) produces an infinite family of $(l+1)$-regular Ramanujan graphs. However, we shall see that the constraints imposed on the parameters $\{p,l\}$ by the LPS and Pizer constructions do {\em{not}} agree. In Section \ref{subsect:Pizer_primes} we give an explicit comparison between the admissible values of the parameter $p$ in the example when $l=5.$

First we wish to summarize the construction via Pizer \cite{pizer:1998}. In particular we wish to explain the elements of \cite[Theorem 5.1]{pizer:1998}. Details are kept to a minimum; the reader is encouraged to consult {\em{op.cit.}}\ for details, in particular \cite[4.]{pizer:1998}. We mention one feature of Pizer's approach in advance: we shall see that here the graph is given via its adjacency matrix. Note that this is of a different flavor from the LPS case. There the edges of the graph were specified ``locally:'' given a vertex of the graph (as an element of a group in Section \ref{subsubsect:CaylePres} or as a class of lattices in Section \ref{subsubsect:infinitelatticetree}), its neighbors were specified directly. (See Section \ref{app:sect:explicitisom} for an explicit parametrization of the edges at a vertex.) In Pizer's approach the adjacency matrix, a Brandt matrix (associated to an Eichler order in the quaternion algebra) specifies the edge structure of the graph. 

\subsection{Overview of the construction}

Let us fix $B=B_{p,\infty}$ to be the quaternion algebra over $\QQ$ that is ramified precisely at $p$ and at infinity. We shall consider orders $\O$ of level $N=pM$ and $N=p^2M$ in $B,$ where $M$ is coprime to $p$.
The vertex set of our graph $G(N,l)$ shall be in bijection with (a subset of) the class group of $\O .$ The class number of $\O$ depends only on the level of the order and hence we may write $H(pM)$ or $H(p^2M)$ for the size of such a graph. In the case where $M=1$ by the Eichler class number formula \cite[Proposition 4.4]{pizer:1998} we have:
\begin{equation}\label{eq:clnum_Hp}
H(p)=\frac{p-1}{12}+\frac{1}{4}\left(1-\left(\frac{-4}{p}\right)\right)+\frac{1}{3}\left(1-\left(\frac{-3}{p}\right)\right);
\end{equation}
\begin{equation}\label{eq:clnum_Hpsq}
H(p^2)=\frac{p^2-1}{12}+\left\lbrace \begin{array}{ll}
0 & \text{ if }p\geq 5\\
\frac{4}{3} & \text{ if }p=3
\end{array}\right.
\end{equation}
where $\left(\frac{\cdot}{\cdot}\right)$ is the Kronecker symbol.

The vertex set of $G(N,l)$ shall have $H(N)$ elements when $N=pM$ and when $N=p^2M$ and $l$ is a quadratic nonresidue modulo $p.$ (Note that in this case the graph $G(p^2M,l)$ is bipartite.) For $N=p^2M$ and $l$ a quadratic residue modulo $p$ the graph $G(p^2M,l)$ is non-bipartite of size $\frac{H(p^2M)}{2}.$ Recall that a similar dichotomy (between bipartite and non-bipartite cases) exists in the LPS construction as well. The following table summarizes the size of $G(p,l)$ and $G(p^2,l)$ for the case where $\left(\frac{l}{p}\right)=1$ (and $p>3$).
\begin{equation}\label{eq:classnumtable}
\begin{array}{c|c|c}
p\mod{12} & H(p) & \frac{H(p^2)}{2} \\ \hline
1 & \frac{p-1}{12} & \multirow{4}{*}{ $\smash[c]{\frac{p^2-1}{12}}$ } \\ 
5 & \frac{p+7}{12} & \\ 
7 & \frac{p+5}{12} & \\ 
11 & \frac{p+13}{12} & \\ 
\end{array}
\end{equation}
The edge structure of the graph $G(N,l)$ is determined via the adjacency matrix. Recall that the rows and columns of the adjacency matrix of a graph are indexed by the vertex set. One entry of the matrix determines the number of edges between the vertices corresponding to its indices. The edge structure of $G(N,l)$ is given by a {\em{Brandt matrix}}. There is a space of modular forms associated to the order $\O$ of the quaternion algebra. This space has dimension as in \eqref{eq:classnumtable} and it carries the action of a Hecke algebra. For every integer $l$ (coprime to $p$) the Brandt matrix $B(N,l)$ describes the explicit action of a particular Hecke operator ($T_l$) on this space.

Restrictions on the parameters $p$ and $l$ guarantee that $B(N,l)$ is in fact the adjacency matrix of a graph. Properties of the resulting graph (e.g.\ the graph being simple and connected, as well as statements about its spectrum and girth) can be phrased as statements about the Brandt matrices $B(N,l)$ and in turn studied as statements about modular forms.

To ensure the edges of the graph $G(N,l)$ are undirected, $B(N,l)$ must be symmetric. By \cite[Proposition 4.6]{pizer:1998} this is the case for $N=pM$ if $p\equiv 1\mod{12}$ and for $N=p^2M$ if $p>3.$

To ensure the graph has no loops we must have $\tr B(N,l)=0,$ and for no multiple edges $\tr (B(N,l))^2=0.$ By \cite[Proposition 4.8]{pizer:1998} these translate to the conditions $\tr B(N,l)=0,$ $\tr B(N,l^2)=H(N).$ (This depends on the relationship of the traces within a family of Brandt matrices $B(N,l)$ for fixed $N$ and varying $l.$) These traces can be given in terms of 
parameters dependent on the order $\O$ \cite[Proposition 4.9]{pizer:1998}.

It turns out that the above conditions together already guarantee that $B(N,l)$ determines a Ramanujan graph. This is the content of the following theorem.

\begin{theorem}\label{thm:Pizer_main}
\cite[Theorem 5.1]{pizer:1998} Let $l$ be a prime coprime to $pM$ and let $N=pM.$ Consider the graph $G(N,l)$ determined by the Brandt matrix $B(N,l)$ as its adjacency matrix. Assume that $B(N,l)$ is symmetric, $\tr B(N,l)=0$ and $\tr B(N,l^2)=H(N).$ Then $G(N,l)$ is a non-bipartite $(l+1)$-regular simple Ramanujan graph on $H(N)$ vertices.

Similarly, let $N=p^2M$ and assume the above conditions $\tr B(N,l)=0$ and $\tr B(N,l^2)=H(N)$ hold. If $l$ is a quadratic nonresidue modulo $p$ then $B(N,l)$ is the adjacency matrix of a bipartite $(l+1)$-regular simple Ramanujan graph on $H(N)$ vertices. If $l$ is a quadratic residue modulo $p$ then $B(N,l)$ is the adjacency matrix of two copies of an $(l+1)$-regular simple non-bipartite Ramanujan graph on $\frac{H(N)}{2}$ vertices.
\end{theorem}

Recall that the quaternion algebra $B$ underlying the construction above is ramified at exactly two places, $p$ and $\infty .$ This uniquely determines the algebra $B=B_{p,\infty }$ (cf.\ \cite[Proposition 4.1]{pizer:1998}). Given a specific $l$ one may ask for what $p$ primes and $N=p$ are the conditions $\tr B(N,l)=0$ and $\tr B(N,l^2)=H(N)$ satisfied. This can be answered by translating the conditions to modular conditions on $p.$ This is carried out for $l=2$ in \cite[Example 2]{pizer:1998}. In the LPS construction above we were interested in $l+1$ regular graphs where $l\equiv 1\mod{4}.$ To compare the families of Ramanujan graphs emerging from the two constructions, in the next section we carry out the same computation for $l=5.$

\subsection{The size of a six-regular Pizer graph}\label{subsect:Pizer_primes}

We wish to consider a special case of Pizer's construction in \cite[Section 5]{pizer:1998}  where the order $\O$ is a (level $p$) maximal order in $B_{p,\infty }$ and the Ramanujan graph is $l+1$ regular. In particular, we are interested in the case where $l=5.$ (Since the LPS construction discussed in Section \ref{subsect:LPS} requires $l\equiv 1\mod{4},$ this is the smallest $l$ where a comparison can be made.) In this section we follow the methods of \cite[Example 2]{pizer:1998} to give explicit modular conditions on $p$ to satisfy Pizer's construction. The Brandt matrix $B(p;5)$ associated to the maximal order $\O\subset B_{p,\infty }$ (of level $p$) is a square matrix of size $H(p).$ It follows from Theorem \ref{thm:Pizer_main} \cite[Proposition 5.1]{pizer:1998} that it is the adjacency matrix of a $6$-regular simple Ramanujan graph if the following conditions hold:
\begin{enumerate}
\item \label{cond:mod12_symmetric} $p\equiv 1\ \mod{12}$
\item \label{cond:tracefirst} $\tr B(p,5)=0$
\item \label{cond:tracesecond} $\tr B(p,5^2)=\Cl \O $
\end{enumerate}
Note that here Condition \ref{cond:mod12_symmetric} guarantees that the graph is symmetric, and Condition \ref{cond:tracefirst} that it has no loops. By \cite[Proposition 4.4]{pizer:1998} the condition $p\equiv 1\ \mod{12}$ gives $\Cl(\O)=Mass \O =\frac{p-1}{12}.$

The Conditions \ref{cond:tracefirst} and \ref{cond:tracesecond} concern the trace of the Brandt matrices $B(p,5)$ and $B(p,25)$ associated to $\O$ of level $p.$ These can be computed using \cite[Proposition 4.9]{pizer:1998}. In particular, {\em{loc.\ cit.}}\ guarantees that Conditions \ref{cond:tracefirst} and \ref{cond:tracesecond} hold under certain conditions. To state these conditions we must introduce some notation. For $m=5$ and $m=25$ respectively, let $s$ be an integer such that $\Delta =s^2-4m$ is negative. Let $t$ and $r$ be chosen such that
\begin{equation}
\Delta =s^2-4\cdot m=\left\lbrace \begin{array}{ll}
t^2r & 0>r\equiv 1 \ \mod \ 4\\
t^24r & 0>r\equiv 2,3 \ \mod \ 4
\end{array}\right.
\end{equation}
Let $f$ be any positive divisor of $t$ and $d:=\frac{\Delta}{f^2}.$ Let $c(s,f,p)$ denote the number of embeddings of $\O_p^d$ into $\O_p$ that are inequivalent modulo the unit group $U(\O_p).$
By \cite[Proposition 4.9]{pizer:1998} we have that
\begin{equation}\label{eq:firstcond_misp}
\text{Condition \ref{cond:tracefirst} is satisfied} \Longleftrightarrow c(s,f,p)=0 \text{ for every $s,f$ with $m=5$ }
\end{equation}
\begin{equation}\label{eq:secondcond_mispp}
\text{Condition \ref{cond:tracesecond} is satisfied} \Longleftrightarrow c(s,f,p)=0 \text{ for every $s,f$ with $m=5^2$}
\end{equation}
The integers $c(s,f,p)$ are given in tables in \cite[pp. 692-693]{pizer1976representability}. We use information in these tables to translate the conditions \eqref{eq:firstcond_misp} and \eqref{eq:secondcond_mispp} into modular conditions on $p.$

First, if $m=5$ the possible values of $s,\Delta ,r,t$ and $f$ are as follows:
$$\begin{array}{c||c|c|c|c|c|c}
s & 0& 1& \multicolumn{2}{|c|}{2}& 3& 4 \\ \hline
\Delta & -20 & -19 & \multicolumn{2}{|c|}{-16}& -11& -4 \\ \hline
t & 1& 1&  \multicolumn{2}{|c|}{2}& 1& 1\\ \hline
r & -5 & -19 & \multicolumn{2}{|c|}{-1} & -11 & -1 \\ \hline
f & 1 & 1& 1 & 2 & 1 & 1 \\ \hline
d & -20& -19& -16&  -4& -11& -4
\end{array}$$
It follows from Condition \ref{cond:mod12_symmetric} that $p\nmid d =\frac{\Delta }{f^2}.$ It follows from the tables in \cite[pp.\ 692--693]{pizer1976representability} that $c(s,f,p)=c(s,f,p)_{p^{2\cdot 0+1}}=0$ if and only if  $d$ is the square of a unit in $\ZZ_p,$ i.e.\ a quadratic residue modulo $p.$ By Condition \ref{cond:mod12_symmetric} we certainly have $\left( \frac{-4}{p}\right)= \left( \frac{-16}{p}\right)=1$ and by quadratic reciprocity $\left( \frac{d}{p}\right)=1$ is equivalent to $\left( \frac{p}{d}\right)=1.$ It follows that by \eqref{eq:firstcond_misp} that Condition \ref{cond:tracefirst} is satisfied if in addition to Condition \ref{cond:mod12_symmetric} $p$ satisfies the following modular conditions.
\begin{equation}\label{eq:modconds_firstpower}
\begin{array}{c|c||l}
c(s,f,p) & \Delta = d & \text{condition}\\ \hline \hline
c(0,1,p)& -20 & p\in \{1,4\}\mod{5} \\ \hline
c(1,1,p)& -19 & p\in \{1,4,5,6,7,9,11,16,17\}\mod{19} \\ \hline
c(3,1,p)& -11 & p\in \{1, 3, 4, 5, 9\} \mod{11}
\end{array}
\end{equation}

Second, to guarantee that the conditions in \eqref{eq:secondcond_mispp} are satisfied, let $m=25.$ Then the possible values of $s,\Delta ,r,t$ and $f$ are as follows:
\begin{equation}\label{eq:sfdelta_squarecond}
\begin{array}{c||c|c|c|c|c|c|c|c|c|c}
s & 0 & 1 & 2& 3 & 4 & 5 & 6 & 7 & 8 & 9 \\ \hline
\Delta  & -100 & -99 & -96 & -91 & -84 & -75 & -64 & -51 & -36 & -19 \\ \hline
t & 5 & 3 & 4 & 1 & 1 & 5 & 4 & 1 & 3 & 1 \\ \hline
r & -1 & -11 & -6& -91 & -21 & -3 & -1 & -51 & -1 & -19 \\ \hline
f & 1,5 & 1,3 & 1,2,4 & 1 & 1 & 1,5 & 1,2,4 & 1 & 1,3 & 1
\end{array}
\end{equation}
By \eqref{cond:mod12_symmetric} and \eqref{eq:modconds_firstpower} we have that $p\nmid d$ for any of the above values of $\Delta $ and $d=\frac{\Delta }{f^2}.$ Then it again follows from the tables in \cite[pp.\ 692--693]{pizer1976representability} that \eqref{eq:secondcond_mispp} is satisfied if and only if for any such $d$ $\left( \frac{d}{p}\right)=1$ or, equivalently by \eqref{cond:mod12_symmetric}, $\left( \frac{p}{d}\right)=1.$ By properties of the Legendre symbol and the previously imposed conditions on the residue class of $p$ modulo $12,$ $5,$ $11$ and $19$ this is true for $\Delta \in \{-100,-99,-75, -64, -36, -19\}$. The remaining cases amount to the following additional modular conditions on $p:$
\begin{equation}\label{eq:modconds_secondpower}
\begin{array}{c|c||l}
\Delta & d=\frac{\Delta }{f^2}& \text{condition}\\ \hline \hline
-96 & -96,\ -24\text{ or }-6& p\in \{1,7\} \mod{8}\\ \hline
-51 & -51=-3\cdot 17 & p\in \{1,2,4,8,9,13,15,16\} \mod{17} \\ \hline
-84 & -84=-12\cdot 7 & p\in \{1,2,4\} \mod{7}\\ \hline
-91 & -91=-7\cdot 13 & p\in \{1,3,4,9,10,12\} \mod{13}
\end{array}
\end{equation}

We summarize the modular conditions on $p$ in the following corollary.
\begin{corollary}\label{cor:6regPizer}
The Brandt matrix $B(p;5)$ associated to a maximal order in $B_{p,\infty }$ by Pizer \cite{pizer:1998} is the adjacency matrix of a $6$-regular simple, connected, non-bipartite Ramanujan graph if and only if $p$ satisfies the following congruence conditions:
\begin{equation}
\begin{array}{c|l}
\text{Modulus} & \text{Remainders allowed} \\ \hline \hline
24 & 1\\ \hline
5 & 1,4\\ \hline
7 & 1, 2, 4\\ \hline
11 & 1, 3, 4, 5, 9\\ \hline
13 & 1, 3, 4, 9, 10, 12\\ \hline
17 & 1, 2, 4, 8, 9, 13, 15, 16\\ \hline
19 & 1, 4, 5, 6, 7, 9, 11, 16, 17\\
\end{array}
\end{equation}
These conditions are equivalent to saying that $p\equiv 1 \mod{24}$ and $p$ is a quadratic residue modulo the primes $5, 7, 11, 13, 17, 19.$ Note that $p$ may belong to one of $1\cdot 2\cdot 3\cdot 5\cdot 6\cdot 8\cdot 9=12\ 960$ residue classes modulo $24\cdot 5\cdot 7\cdot  11\cdot 13\cdot 17\cdot 19=38\ 798\ 760.$
\end{corollary}

The Corollary describes the set of primes $p$ for which $G(p,5)$ is a six-regular Ramanujan graph. The condition $p\equiv 1\mod{4},$ $p\equiv 1,4\mod{5=l}$ guarantees that for these primes the LPS construction is a six-regular graph as well.

\begin{remark}\label{rmk:primeswhenbothexist}
The smallest prime satisfying all the congruence conditions of Corollary \ref{cor:6regPizer} is $53881.$ This corresponds to a $6$-regular Pizer graph with $4490$ vertices. Amongst the first one million primes, $1670$ satisfy all these congruence conditions.
\end{remark}

\section{Relationship between LPS and Pizer constructions}\label{subsect:LPS_Pizer_compare}

We wish to compare the two different approaches to constructing Ramanujan graphs that we have discussed. Throughout the previous sections, we have seen that the constructions of LPS and Pizer (recall the latter agree with supersingular isogeny graphs for particular choices) have similar elements. In this section, we wish to further highlight these similarities, as well as the discrepancies between the two approaches.

First let us revisit the chains of graph isomorphisms/bijections that the respective constructions fit into. These are as follows:
\footnotesize
\begin{equation*}\label{eq:chains_compare}
\begin{split}
\text{\bf \normalsize (LPS) } \Cay(\PSLt(\FF_p),S)\cong \Gamma(2p)\backslash \PGLt(\Ql)/\PGLt(\Zl)\cong & G'(\QQ)\backslash H_{2p}(\A_f)/K_0^{2p} \\
(\O[l^{-1}])^{\times }\backslash \GLt(\Ql)/\GLt(\Zl) \cong & B^{\times }(\QQ)\backslash B^{\times }(\A_f)/B^{\times }(\hat{\ZZ})\cong \Cl \O \cong \text{SSIG } \text{\bf \normalsize (Pizer)}
\end{split}
\end{equation*}
\normalsize
Recall that in the first line, we have the LPS construction in terms of a Cayley graph on the group $\PSLt(\FF_p);$ it corresponds to the ``local double coset graph'' defined by taking a finite quotient of an infinite tree of homothety classes of lattices. The vertex set of this graph is in bijection with the adelic double cosets on the right-hand side. (For the sake of this comparison we omitted the infinite place.)

On the right-hand end of the second line, we have the supersingular isogeny graphs discussed in Part \ref{part1}. These are symmetric simple graphs isomorphic to $G(p,l)$ constructed by Pizer (see Section \ref{sect:Pizer}) when $p\equiv 1\mod{12}.$ The vertex set of $G(p,l)$ is the class group of a maximal order $\O$ in the quaternion algebra $B_{p,\infty }.$ This set is in bijection with the adelic double cosets. Via strong approximation (see Section \ref{subsubsect:SA_Pizer}) these adelic double cosets are in bijection with local double cosets, which at a place $l$ where $B_{p,\infty }$ splits can be written as the left-hand side object.

Despite the similarities between these chains of bijections, there are significant discrepancies between the two objects. First of all, there is a discrepancy in the underlying quaternion algebras. For the LPS graphs we considered the underlying algebra of Hamiltonian quaternions ($B_{2,\infty }$). Varying the parameter $p$  we get different Ramanujan graphs by changing the congruence subgroup $\Gamma(2p)$  without ever changing the underlying algebra. On the other hand the Pizer graphs were constructed using $B=B_{p,\infty }.$ The underlying quaternion algebra varies with the choice of the parameter $p.$ We note that  the construction in LPS can be carried out  for any $B$ ramified at $\infty$ and split at $l$, and would still result in Ramanujan graphs  (see \cite[Theorem 7.3.12]{lubotzky:2010}). However, in this more general case we do not have a clear path for obtaining an explicit description of these graphs as Cayley graphs. For additional details see \cite[Remark 7.4.4(iv)]{lubotzky:2010}. If one took $B_{p,\infty}$ for both the LPS and Pizer cases, the infinite families of Ramanujan graphs formed would differ because the LPS family is formed by varying the subgroup  $\Gamma(2p)$ (or more generally $\Gamma(N)$ for $l$ a quadratic residue mod $N$) while the Pizer family is formed by varying the quaternion algebra $B_{p,\infty}$.

Let us consider the choice of parameters next. For the LPS graphs we required only that $l\equiv 1\mod{4}$ and that $p$ is odd and prime to $l.$ If $-1$ is a quadratic residue modulo $p$ then the resulting graph is isomorphic to a subgroup of $\PGLt(\ZZ/p\ZZ)$ \cite[Theorem 7.4.3]{lubotzky:2010}. Furthermore, if $l$ is a quadratic residue modulo $2p$ then this graph is non-bipartite and isomorphic to the Cayley graph of $\PSLt(\FF_p)$ with $\frac{p^3-p}{2}$ elements.

In the case of the Pizer graphs $G(N,l)$ we must have $N=pM$ coprime to $l.$ Further congruence conditions on $N$ guarantee properties of the resulting graph (see Section \ref{sect:Pizer}), e.g.\ $p\equiv 1\mod{12}$ guarantees that the adjacency matrix is symmetric. The number of vertices in $G(N,l)$ is then $H(N),$ the class number of an order of level $N$ in $B_{p,\infty }.$ For example if $N=p\equiv 1\mod{12},$ then this results in a graph of size $\frac{p-1}{2}.$

To compare the two in the simplest case when $l\equiv 1\mod{4}$, i.e.\ $l=5,$ recall that Corollary \ref{cor:6regPizer} gives the exact congruence conditions on $p$ so that the Pizer construction of the graph $G(p,5)$ is a six-regular Ramanujan graph on $\frac{p-1}{12}$ vertices. For these primes, the LPS construction also produces a Ramanujan graph. The size of the two graphs is very different. Notice however that when both graphs exist the size of the LPS graph is divisible by the size of the Pizer graph (cf. Remark \ref{rmk:primeswhenbothexist}).

Let us turn our attention to the local double coset objects in the above chain of bijections. In the second line, corresponding to Pizer graphs, we have $(\O[l^{-1}])^{\times }$ appearing where $\O $ is an order of the quaternion algebra $B_{p,\infty }.$ For the graph $G(p,l)$ this $\O $ is an order of level $p,$ i.e.\ a maximal order. The corresponding subgroup $(\O[l^{-1}])^{\times }$ of $B^{\times }(\ZZ[l^{-1}])$ is analogous to the subgroup $\Gamma =G'(\ZZ[l^{-1}])$ for the LPS construction. This is much larger than the congruence subgroup $\Gamma (2p)\leq \Gamma $ that appears in the local double coset objects in that case.

The fact that the LPS construction involves this {\em{smaller}} congruence subgroup $\Gamma (2p)$ also accounts for the discrepancy between the two lines at the adelic double cosets. Recall from Section \ref{subsubsect:SA_LPS} that $H_{2p}$ was not the entire $G'(\A)$ but instead a finite index normal subgroup of it. We note that if one replaced $\Gamma(2p)$ in the LPS construction with $\Gamma(2N)$, where $p \mid N$, the  LPS graph $\Gamma(2N)\backslash \PGLt(\Ql)/\PGLt(\Zl)$ is a finite cover of $\Gamma(2p)\backslash \PGLt(\Ql)/\PGLt(\Zl)$ \cite[Section 3]{li:1996}.

One may wonder if an object analogous to $\Cl(\O)$ could be appended to the chain of bijections for LPS graphs. Or even if, in the local double coset object for LPS graphs $\Gamma (2p)$ could be written as $(\O_{2p}(\ZZ[l^{-1}]))^{\times }$ as well, for a quaternion order $\O_{2p}.$ (More precisely, if $\Gamma (2p)$ agrees with the image of $(\O_{2p}(\ZZ[l^{-1}]))^{\times }$ under the map $B^{\times }\rightarrow G'$ for some order $\O_{2p}.$)

The answer to the second question is affirmative. Using the basis $1,{\mathbf{i}}, {\mathbf{j}},{\mathbf{k}}$ for $B=B_{2,\infty }$ the requisite relationship holds between $\O_{2p}$ and $\Gamma (2p)$ for the order $\O_{2p}$ spanned by $\{1,2p\vi, 2p\vj,2p\vk \}.$ Note that this order has level $2^5p^3,$ hence it is not an Eichler order.

We remark that the size of the class set of this $\O_{2p}$ can be computed using \cite[Theorem 1.12]{pizer1980algorithm} and it turns out to be $\frac{4p^2(p+1)+4}{3}$ or $\frac{4p^2(p+1)}{3}$ if $p\equiv 1\mod{3}$ or $p\equiv 2\mod{3}$ respectively. This is clearly different from the size of $\PSLt(\FF_p)$ which is  a numerical obstruction to extending the chain of isomorphisms for LPS graphs analogously to the row for Pizer graphs.

\bigskip

\providecommand{\bysame}{\leavevmode\hbox to3em{\hrulefill}\thinspace}
\providecommand{\MR}{\relax\ifhmode\unskip\space\fi MR }
\providecommand{\MRhref}[2]{%
  \href{http://www.ams.org/mathscinet-getitem?mr=#1}{#2}
}
\providecommand{\href}[2]{#2}

\end{document}